\documentclass[12pt,a4paper]{article}

\usepackage{amsmath} 
\usepackage{latexsym} 
\usepackage{amssymb}
\usepackage{amsthm}

\DeclareMathOperator*{\esssup}{ess\,sup}

\newcommand {\eqdef} {\stackrel{\rm def}{=}} 

\newtheorem{theorem}{Theorem}[section]
\newtheorem{lemma}[theorem]{Lemma}
\newtheorem{corollary}[theorem]{Corollary}
\newtheorem{definition}[theorem]{Definition}
\newtheorem{example}[theorem]{Example}
\newtheorem{remark}[theorem]{Remark}

\usepackage{hyperref}

\begin{document}

\title{
Conditions for correct solvability of a first order linear differential equation in space $L_p(R)$
and asymptotic properties of its solutions}

\author{
M. Lukachev
\and
L. Shuster\thanks{L. Shuster, Department of Mathematics, Bar-Ilan University, 52900 Ramat Gan, Israel;}
}

\date{}
\maketitle

\begin{abstract}
In this work we consider the differential equation
\begin{equation}\label{abstr:theequation}
-y'(x)+q(x)y(x) = f(x), \quad x\in R,
\end{equation}
where $f\in L_p(R),\ p\in[1,\infty] \ (L_\infty(R):=C(R))$
and $0\le q\in L^{loc}_1(R)$.
The equation \eqref{abstr:theequation} is called correctly solvable in the given space $L_p(R)$ if for any $f\in L_p(R)$ there is a unique solution $y\in L_p(R)$ and the following inequality
\begin{equation*}
\|y\|_p \le c(p)\|f\|_p, \quad \textrm{for all} \quad f\in L_p(R),
\end{equation*}
holds with absolute constant $c(p) \in (0,\infty)$. We find a criterion for correct solvability of the above equation in space $L_p(R)$ and study the asymptotic properties of its solutions.
\\ \\ \\ \noindent {\bf MSC~2000}: 34B05, 34B40, 34C11.
\vspace{1mm}
\\ \noindent {\bf Key Words}: first order linear differential equation, correct solvability, asymptotic properties of solutions.
\end{abstract}

\newpage

\section{Introduction}
\renewcommand{\theequation}{\arabic{section}.\arabic{equation}}
\setcounter{equation}{0}

\subsection{Basic definitions and preliminaries}\label{Preliminaries} 

Below we consider the differential equation
\begin{equation}\label{theequation}
-y'(x)+q(x)y(x) = f(x), \quad x\in R,
\end{equation}
where
\begin{equation*}
f\in L_p(R),\ p\in[1,\infty] \quad (L_\infty(R):=C(R)),
\end{equation*}
and
\begin{equation}\label{q} 
0\le q\in L^{loc}_1(R).
\end{equation}
 
As a solution of equation \eqref{theequation} we assume any absolutely
continuous function $y$ that satisfies \eqref{theequation} almost everywhere on $R$.
We need the following definition:
\begin{definition}
We call the equation \eqref{theequation} correctly solvable in the given space $L_p(R)$ if the following two assertions hold (\cite[ch.III, 6.2]{Co}):
\begin{enumerate}
\item[I)] for any $f\in L_p(R)$, the equation \eqref{theequation} has a unique solution
$y\in L_p(R)$;
\item[II)] the solution $y\in L_p(R)$ of equation \eqref{theequation} satisfies the inequality
\begin{equation*}
\|y\|_p \le c(p)\|f\|_p,\quad \textrm{for all }\ f\in L_p(R),
\end{equation*}
\end{enumerate}
where $c(p)$ is an absolute positive constant.
\end{definition}
We note that instead of ``correct solvability'' other equivalent terms also are used (see \cite[ch.5, \S50-51]{MS}).

Throughout the paper in accordance with I)-II), we denote by the symbol $y$ only a solution of equation \eqref{theequation} that belongs to $L_p(R)$.

In what follows, $c,\ c(\cdot),\ c_1,\ c_2,\ \dots$ stand for absolute positive constants, which are not essential for exposition and may differ even within a single chain of calculations.

The equation \eqref{theequation} was considered in papers \cite{CS,CS1,Est} in connection with studying of the singular boundary value problem
\begin{equation}\label{bvp:equation}
-y'(x)+q(x)y(x) = f(x), \quad x\in R,
\end{equation}
\begin{equation}\label{bvp:lim-y}
\lim\limits_{|x|\to \infty} y(x) = 0.
\end{equation}
In particular, the following result was obtained in \cite{CS}:

\begin{theorem}{\rm\cite{CS}}\label{bvp:theorem}
Let $p\in[1,\infty)$. The boundary value problem \eqref{bvp:equation}-\eqref{bvp:lim-y} is correctly solvable in space $L_p(R)$ if and only if there is $a\in (0, \infty)$ such that
\begin{equation}\label{inf_int_q} 
q_0(a) > 0, \quad q_0(a) \eqdef \inf_{x\in R} \int_{x-a}^{x+a} q(t)dt.
\end{equation}
In addition, the problem \eqref{bvp:equation}-\eqref{bvp:lim-y} is correctly solvable in $C(R)$ if and only~if
\begin{equation}\label{lim_int_q}
\lim_{|x|\to \infty} \int_{x-a}^{x+a} q(t)dt = \infty, \quad \textrm{for all }\ a\in (0,\infty).
\end{equation}
\end{theorem}
\noindent
Moreover, as shown in \cite{CS}, a criterion for correct solvability in space $L_p(R),$\ $p\in [1,\infty]$ of the boundary value problem \eqref{bvp:equation}-\eqref{bvp:lim-y} can be stated in terms of some auxiliary function $q^*$ of coefficient $q$. Therefore, let us now introduce this function. For this we temporarily assume that in addition to condition \eqref{q} we have $\|q\|_1 = \infty$. Then for any given $x\in R$, we define a function
\begin{equation}\label{d}
d(x)=\inf_{d>0}\ \Big\{d : \int_{x-d}^{x+d}q(t)dt=2\Big\}.
\end{equation}
It is clear that the function $d(x)$ is defined correctly. (The function  $d(x)$ was first used in \cite{CS1}). Now set $q^{*}(x)=\frac{1}{d(x)}$. Then the following Theorem \ref{bvp:theorem-alt} provides a criterion for correct solvability in space $L_p(R), p\in [1,\infty]$ of the boundary value problem \eqref{bvp:equation}-\eqref{bvp:lim-y} stated in terms of the function $q^{*}$:
\\
\begin{theorem}{\rm \cite{CS}}\label{bvp:theorem-alt} 
Let $p\in[1,\infty)$. The boundary value problem \eqref{bvp:equation}-\eqref{bvp:lim-y} is correctly solvable in space $L_p(R)$ if and only if the following two conditions hold:
\begin{equation}
\int_{-\infty}^{0}q(t)dt = \int_{0}^{\infty}q(t)dt = \infty,
\end{equation}
\begin{equation}
\gamma \eqdef \inf_{x\in R}q^{*}(x) > 0 \quad \big(\textrm{or }\ \sup_{x\in R}d(x) < \infty \big).
\end{equation}
In addition, the problem \eqref{bvp:equation}-\eqref{bvp:lim-y} is correctly solvable in $C(R)$ if and only if
\begin{equation}\label{lim_q*}
q^{*}(x)\to\infty \textrm{ as } |x|\to\infty \quad \big(\textrm{or }\ d(x)\to 0 \textrm{ as } |x|\to\infty \big).
\end{equation}
\end{theorem}

In connection with condition \eqref{lim_int_q} of Theorem \ref{bvp:theorem} and condition \eqref{lim_q*} of Theorem \ref{bvp:theorem-alt} we note the following statement:

\begin{lemma}{\rm [\S\ref{D:proof}]}
The condition \eqref{lim_int_q} and the following condition \eqref{lim_d} are equivalent.
\begin{equation}\label{lim_d}
\lim_{|x| \to \infty} d(x) = 0
\end{equation}
\end{lemma}

Let us now introduce our goals and results.


\subsection{Problem A -- Conditions for correct solvability of a first order linear differential equation}\label{A}

Our main goal is to establish a criterion for correct solvability in space $L_p(R),\ p\in [1, \infty]$ of equation \eqref{theequation} without requirement  \eqref{bvp:lim-y}. We obtained the following result:

\begin{theorem}{\rm [\S\ref{A:proof}]}\label{A:maintheorem}
Let $p\in[1,\infty]$. The equation \eqref{theequation} is correctly solvable in space $L_p(R)$ if and only if condition \eqref{inf_int_q} holds.\\
In particular, only one of the assertions A) and B) holds: \\
A) for all $p \in [1, \infty]$, the equation \eqref{theequation} is correctly solvable in $L_p(R)$; \\
B) for all $p \in [1, \infty]$, the equation \eqref{theequation} is not correctly solvable
in $L_p(R)$.\\
In addition, the solution $y\in L_p(R)$ of equation \eqref{theequation} is of the form:
\begin{equation}\label{y=int-G(x,t)fdt}
y(x) \eqdef (Gf)(x) = \int_{-\infty}^{\infty} G(x,t)f(t)dt,\quad x\in R,\\
\end{equation}
where
\begin{equation}\label{G(x,t)=}
G(x,t) =
\left\{
\begin{array}{l}
0, \quad  t<x,\ x\in R\\
\displaystyle
\exp \Big( -\int_{x}^{t}q(\xi)d\xi \Big),\quad t\ge x,\ x\in R.
\end{array}
\right.
\end{equation}
\end{theorem}
\noindent
In connection with Theorem \ref{A:maintheorem} see \cite{CS-LS}.

Let us briefly comment on the above statement. Theorem \ref{A:maintheorem} shows that in case of equation \eqref{theequation} for the spaces with integral metric in the
classical bundle ``equation - space - boundary conditions'', the last
element is meaningless. Existence, uniqueness, and boundedness of solution
for equation \eqref{theequation} in space $L_p(R),\ p\in[1,\infty]$ is entirely
defined by ``equation - space''. Thereby, the boundary conditions for the correctly solvable equation \eqref{theequation} in space $L_p(R),\ p\in[1,\infty)$ are not of importance because according to Theorem \ref{bvp:theorem} they hold automatically, i.e., with no additional assumptions. (It is quite possible that this is the reason why in the problem \eqref{bvp:equation}-\eqref{bvp:lim-y} the number of boundary conditions is greater than the order of equation.) From the given theorems, it also follows that for equation \eqref{theequation} a transition from integral metric to the uniform metric makes the boundary value conditions of significant importance: the solution of correctly solvable equation \eqref{theequation} in the space $C(R)$ satisfies \eqref{bvp:lim-y} regardless of $f\in C(R)$ only if \eqref{lim_int_q} holds, which is stronger than
\eqref{inf_int_q}.

\begin{corollary}
Suppose that condition \eqref{inf_int_q} holds. Then solution $y\in L_p(R)$, $p\in[1,\infty)$ of equation \eqref{theequation} satisfies the inequality
\begin{equation}\label{A:cor1}
\| q^{1/p}y\|_p\ \le\ c \|f\|_p,\quad \it{for\ all\ }\ f \in L_p(R).
\end{equation}
\end{corollary}

\begin{corollary}
Suppose that condition \eqref{inf_int_q} holds. Then equation \eqref{theequation} is separable in the space  $L_1(R)$, i.e., the solution $y\in L_1(R)$ of equation \eqref{theequation} satisfies the inequality
\begin{equation}\label{A:cor2}
\|y'\|_1 + \|qy\|_1\ \le\ 3\|f\|_1,\quad \it{for\ all\ }\ f\in {L_1}(R).
\end{equation}
\end{corollary}

\begin{remark}
{\rm The way we obtain \eqref{A:cor1} and \eqref{A:cor2} is the same as in \cite{CS} and therefore, is not given here. The problem of separability of  differential equations was first studied in \cite{EG1,EG2}. For detailed analysis of condition \eqref{inf_int_q} and examples, see \cite{CS}.}
\end{remark}

\begin{example}
Let
\begin{equation*}
q(x) = e^{x^2} + e^{x^2} \cos{e^{x^2}},\quad x\in R.
\end{equation*}
Then equation \eqref{theequation} is correctly solvable in space $L_p(R),\ p\in [1,\infty]$.
\end{example}
\noindent
Now we would like to study the properties of solutions $y\in L_p(R), p\in [1,\infty]$ of equation \eqref{theequation}, especially their behavior in infinity. 
By difficulty of integral representation of solution $y$, its asymptotic behavior is of a practical interest.

\begin{remark}
For all problems stated below we assume permanently that condition \eqref{inf_int_q} holds, i.e. the equation \eqref{theequation} is correctly solvable in the given space $L_p(R),\ p\in [1,\infty]$.
\end{remark}


\subsection{Problem B -- The asymptotic majorant of solutions for a linear differential equation of a first order}\label{B}

For our next problem, we denote by $D_p$ the set of all solutions $y\in L_p(R)$ of equation \eqref{theequation} with right side $f$ which belongs to the unit sphere
$S_p = \{ f\in L_p(R) : \|f\|_p = 1\}$ in space $L_p(R)$, $p \in[1,\infty]$, i.e.
\begin{equation}
D_p = \big\{\ y: -y'(x) + q(x)y(x) = f(x), \ y\in L_p(R),\ \|f\|_p=1\ \big\}.
\end{equation}

We need the following definition:

\begin{definition}\label{definition:majorant}
Let equation \eqref{theequation} be correctly solvable in the given space $L_p(R)$, $p \in[1,\infty]$. A continuous and positive (for all $x\in R$) function $\varkappa_p(x)$ is called the asymptotic majorant for the set of solutions $D_p$ of equation \eqref{theequation}, if the following relation holds:
\begin{equation}
\lim_{|x|\to \infty} \frac{\sup_{y\in D_p^q} |y(x)|}{\varkappa_p(x)} = 1.
\end{equation}
\end{definition}

Our goal here is for fixed $p\in [1, \infty]$ to find an asymptotic majorant $\varkappa_p$ for the set of solutions $D_p$ for correctly solvable equation \eqref{theequation} in space $L_p(R)$ (further we will say shortly ``asymptotic majorant of solutions'').
We note that a problem of asymptotic majorant was studied in \cite{Maj} for the Sturm-Liouville equation. For equation \eqref{theequation} this problem is studied for the first time here. We note that we widely use methods proposed in \cite{Maj}. In particular, as in \cite{Maj}, we primarily find asymptotic estimates for one class of definite integrals (see below \S\ref{B1:proof}-\ref{B2:proof}) and then use known facts about general properties of linear functionals defined in spaces $L_p(R)$ (see \S\ref{B3:proof}).

Let us now note that the function
\begin{equation}\label{G_p}
G_p(x) = \sup_{y\in D_p} |y(x)|
\end{equation}
is, obviously, an asymptotic majorant of solutions $y \in D_p$. Nevertheless, this formula does not solve the problem stated above because it represents the asymptotic majorant of solutions in implicit form. Our goal, however, is to find an explicit form for the function $\varkappa_p(x), p\in [1,\infty]$ from definition \ref{definition:majorant}.

Thus, our result is the following Theorem \ref{B:maintheorem}:

\begin{theorem}{\rm [\S\ref{B3:proof}]}\label{B:maintheorem}
Let condition \eqref{q} for $q$ hold. Suppose that there exist an absolutely continuous function $q_1(x)>0$ and a function $q_2\in L_1^{loc}(R)$ such that
\begin{equation*}
q(x) = q_1(x) + q_2(x),\quad x\in R
\end{equation*}
Let there exist a continuous and positive for all $x\in R$ function $s(x)$ such that the following conditions hold:\\
\begin{enumerate}
\item[a)]
\begin{equation*}
s(x) \to \infty \quad \textrm{as} \quad |x|\to \infty,
\end{equation*}
\item[b)]
\begin{equation*}
\frac{1}{s(x)} \ge \frac{|q'_1(x)|}{q_1^2(x)} \quad \textrm{for all} \quad |x| \gg 1,
\end{equation*}
\item[c)]
\begin{equation*}
\lim_{|x|\to \infty} \frac{s(x)}{xq_1(x)}=0,
\end{equation*}
\item[d)]
for some $\nu\in[1,\infty)$ the inequalities hold:
\begin{equation*}
\frac{1}{\nu}\le\frac{s(t)}{s(x)}\le\nu,\quad t\in \Delta(x),\ |x|\gg 1,
\end{equation*}
where
\begin{equation*}
\Delta(x)=\Big[\Delta^{-}(x),\ \Delta^{+}(x)\Big]=\Big[x-\frac{s(x)}{q_1(x)},\ x+\frac{s(x)}{q_1(x)}\Big],\quad x\in R,
\end{equation*}
\item[e)]
\begin{equation*}
\tilde{\varkappa}(x) \to 0 \ \textrm{ as } \ |x| \to \infty,
\end{equation*}
where
\begin{equation*}
\tilde{\varkappa}(x) = \sup_{t\in \Delta(x)} |\varkappa(t)|, \quad x\in R,
\end{equation*}
\begin{equation*}
\varkappa(t) = q_1(t) \int_{x}^{t} \frac{q_2(\xi)}{q_1(\xi)}d\xi,\quad t\in \Delta(x).
\end{equation*}
\end{enumerate}
Then for $p\in [1, \infty]$ the asymptotic majorant $\varkappa_p(x)$ for solutions $y\in D_p$
\begin{equation*}
D_p = \big\{\ y: -y'(x) + q(x)y(x) = f(x), \ y\in L_p(R),\ \|f\|_p=1\ \big\},
\end{equation*}
is of the form:
\begin{eqnarray*}
\varkappa_p(x)= 
\left\{
\begin{array}{cc}
1, & \quad p=1,\ x\in R, \\
\displaystyle \frac{1}{(p')^{1/p'}} \frac{1}{q_1(x)^{1/p'}}, & \quad
p\in (1,\infty),\ p'=\frac{p}{p-1},\ |x|\gg 1,
\vspace{2mm}\\
\displaystyle \frac{1}{q_1(x)}, & \quad p=\infty, \ |x|\gg 1.
\end{array}
\right.
\end{eqnarray*}
\end{theorem}

\begin{example}
Let
\begin{equation*} 
q(x) = e^{x^2} + e^{x^2} \cos{e^{x^2}},\quad x\in R.
\end{equation*}
Let $\varkappa_p(x)$ be the asymptotic majorant of solutions for equation \eqref{theequation}. Then
\begin{eqnarray*}
\varkappa_p(x) = 
\left\{
\begin{array}{ll}
\quad \quad 1, & \quad p=1,\ x\in R, \\
\displaystyle \frac{1}{(p')^{1/p'}} \frac{1}{e^{x^2/p'}}, & \quad
p\in (1,\infty),\ |x|\gg 1, \\
\displaystyle \quad \ \ \frac{1}{e^{x^2}}, & \quad p=\infty,\ |x|\gg 1.
\end{array}
\right.
\end{eqnarray*}
\end{example}


\subsection{Problem C -- Sharp by order two-sided estimates of asymptotic majorant $\varkappa_p(x)$ of solutions from the set~$D_p$}\label{C}

In Problem B (\S\ref{B}) we have seen that there are quite a lot of requirements for the function $q$, and finding the asymptotic majorant of solutions could be rather difficult task. In those cases we can find estimates for asymptotic majorant itself. Of course, such estimates are less precise, but it is natural to expect that in such case the requirements for the function $q$ will be less strong than in Theorem \ref{B:maintheorem}. We note primarily that from correct solvability in space $L_p(R), p\in [1,\infty]$ of equation \eqref{theequation} it follows that $\|q\|_1=\infty$ (see Theorem \ref{A:maintheorem}). Therefore, the function
\begin{equation*}
d(x)=\inf_{d>0}\ \Big\{d : \int_{x-d}^{x+d}q(t)dt=2\Big\},\quad x\in R
\end{equation*}
is defined correctly (see \eqref{d} in \S\ref{Preliminaries}).

The following Theorem \ref{C:maintheorem} contains sharp by order two-sided estimates of the function $G_p(x)$ (see \eqref{G_p}) in terms of the function $d$.

\begin{theorem}{\rm [\S\ref{C:proof}]}\label{C:maintheorem}
Let there exist $\alpha \ge 1$ and $\beta > 0$ such that for all $|x| \gg 1$ the inequality holds:
\begin{equation}
\frac{1}{\alpha} \le \frac{d(t)}{d(x)} \le \alpha, \quad |t-x| \le \beta,
\end{equation}
Then for $p\in (1,\infty]$ there exist $c,c(p)\in (0, \infty)$ such that for all $x\in R$ the inequalities hold:
\begin{equation}\label{C:result}
c^{-1}d(x)^{1/p'} \le G_p(x) \le c(p)d(x)^{1/p'}.
\end{equation}
\end{theorem}


\subsection{Problem D -- A problem of the $\varepsilon$-strip}\label{D}

We continue to study the asymptotic behavior of solutions $y\in L_p(R),\ p\in[1,\infty]$ of equation \eqref{theequation}. For this purpose let us shortly comment on Theorem \ref{C:maintheorem} (see \S\ref{C}). Suppose that in addition to conditions of this theorem the function $d(x) \to 0$ as $|x|\to \infty$. Then because of \eqref{C:result} the solutions $y\in D_p$ of equation \eqref{theequation} tend to zero as $|x|\to \infty$ uniformly regardless their choice in the set $D_p$. Indeed, the equality
\begin{equation}\label{lim_d(x)}
\lim_{|x|\to \infty} d(x) = 0
\end{equation}
means that for any $\varepsilon>0$ there exists $x_0=x_0(\varepsilon)$ such that for $|x|\ge x_0(\varepsilon)$ the following relations ($c$ is a constant from \eqref{C:result}) hold:
\begin{gather*}
d(x)^{1/p'} \le \frac{\varepsilon}{c}, \\
\sup_{y\in D_p}|y(x)| = G_p(x) \le cd(x)^{1/p'}\le \varepsilon.
\end{gather*}
Thereby, we introduce the following definition:
\begin{definition}\label{definition:tend}
We say that solutions of equation \eqref{theequation} tend in whole to zero as $|x|\to \infty$ if for any $\varepsilon > 0$ there exists $x_0=x_0(\varepsilon) \gg 1$ such that for all $|x|\ge x_0$ the inequality holds:
\begin{equation*}
|y(x)| \le \varepsilon \quad \textrm{ for all }\ y\in D_p,
\end{equation*}
\end{definition}
i.e. for $|x|\ge x_0(\varepsilon)$ all integral curves from $D_p$ are contained in strip
\begin{equation*}
B = \big\{\ (x,y): |y|\le \varepsilon,\ |x|\ge x_0(\varepsilon)\ \big\} 
\end{equation*}
of the plane $XOY$.

Thus, our very next goal is to find a complete answer for the following question: is the stated in definition \ref{definition:tend} property of solutions $y\in D_p$ follows from both condition \eqref{C:result} of Theorem \ref{C:maintheorem} and condition \eqref{lim_d(x)} or, perhaps, this property is entirely defined by requirement \eqref{lim_d(x)} only? The solution of this question is provided by the following Theorem \ref{D:maintheorem}:

\begin{theorem}{\rm [\S\ref{D:proof}]}\label{D:maintheorem}
For $p=1$ the solutions of equation \eqref{theequation} do not tend in whole to zero as
$|x| \to \infty$. For $p\in (1, \infty]$ the solutions of equation \eqref{theequation} tend in whole to zero as $|x| \to \infty$ if and only if one of the equivalent conditions \eqref{lim_int_q} or \eqref{lim_d} holds.
\end{theorem}


\subsection{Problem E -- Criterion of compactness in space $L_p(R)$ for resolvent of differential operator of a first order}\label{E}

To present this problem we give some definitions. Fix $p\in [1, \infty]$. Let $H_p$ be the following set:
\begin{equation*}
H_p = \left\{\ y\in L_p : -y'(x) + q(x)y(x) = f(x), \quad \forall f\in L_p \ \right\},
\end{equation*}
i.e. a set comprised of solutions $y\in L_p(R)$ of equation  \eqref{theequation} with right side
$f\in L_p(R)$.
We remind that throughout below the equation \eqref{theequation} assumed to be correctly solvable in space $L_p(R)$ and do not mention this requirement anymore. 
Let us denote $L$ be a differential operator given on the set $H_p$:
\begin{equation*}
L = -\frac{d}{dx} + q(x),\quad x\in R.
\end{equation*}
Since equation \eqref{theequation} is correctly solvable, for operator $L$ there exists an inverse operator $L^{-1}$. Thereby (see a definition of correct solvability)
\begin{enumerate}
\item[1)] Operator $L^{-1}$ is defined on entire space $L_p(R)$;
\item[2)] Operator $L^{-1}: L_p(R) \to L_p(R)$ is bounded;
\item[3)] Operator $L^{-1}$ operates on $f\in L_p(R)$ by the following rule:
\begin{equation*}
(L^{-1}f)(x):=(Gf)(x) \eqdef \int_{x}^{\infty}e^{-\int\limits_{x}^{t}q(\xi)d\xi}f(t)dt, \quad x\in R.
\end{equation*}
\end{enumerate}
The operator $L^{-1}(=G)$ is the resolvent of operator $L$ and was first considered in \cite{Comp} in connection with condition
\begin{equation}\label{E:qge1}
1 \le q\in L_1^{loc}(R).
\end{equation}
In particular, in \cite{Comp} the following result was obtained:
\begin{theorem}{\rm \cite{Comp}}\label{E:prevtheorem}
Let condition \eqref{E:qge1} hold. For fixed $p\in [1, \infty]$ the operator $L^{-1} : L_p(R) \to L_p(R)$ is compact if and only if condition \eqref{lim_int_q} holds.
\end{theorem}

\begin{remark}
{\rm The requirement \eqref{lim_int_q} is called a condition of A.M.Molchanov (see \cite{SH-M,N}).}
\end{remark}

Our goal here is to get a criterion of compactness of operator $L^{-1}$ in our case, for apriori condition \eqref{q}. In addition, it would be interesting to know how this condition is connected to other problems that we solved.
Our result is the following Theorem \ref{E:maintheorem}:
\begin{theorem}{\rm [\S\ref{E:proof}]}\label{E:maintheorem}
Let condition \eqref{q} hold. Then for fixed $p\in[1,\infty]$ the operator $L^{-1}:L_p(R)\to L_p(R)$ is compact if and only if condition \eqref{lim_int_q} (or equivalent condition \eqref{lim_d}) holds.
\end{theorem}

Theorem \ref{E:maintheorem} shows that for correctly solvable equation \eqref{theequation} in space $L_p(R), p\in (1,\infty]$ the problems D and E are equivalent.


\section{Solution of Problem A}\label{A:proof}
\setcounter{equation}{0}

In this section we prove Theorem \ref{A:maintheorem} which is the main result of our work (see \S\ref{A}). For convenience we remind below its statement. 
\\ \\
{\bf Theorem \ref{A:maintheorem}.}
{\it
Let $p\in[1,\infty]$. The equation \eqref{theequation} is correctly solvable in space $L_p(R)$ if and only if there is $a\in (0, \infty)$ such that
\begin{equation}\label{A:inf_int_q}
q_0(a) > 0, \quad q_0(a) \eqdef \inf_{x\in R} \int_{x-a}^{x+a} q(t)dt.
\end{equation}
In particular, only one of the assertions A) and B) holds: \\
A) for all $p \in [1, \infty]$, the equation \eqref{theequation} is correctly solvable in space $L_p(R)$; \\
B) for all $p \in [1, \infty]$, the equation \eqref{theequation} is not correctly solvable
in space $L_p(R)$.
}\\

\subsection{Proof of main result for $p\in [1,\infty)$.}\label{A:case_p<infty} 

\paragraph{Proof of Theorem \ref{A:maintheorem}. \it{Necessity.}}
Below in the ``necessity part'' of Theorem \ref{A:maintheorem}, we assume that equation \eqref{theequation} is correctly solvable in space $L_p(R)$ for given $p\in [1,\infty)$. We need the following lemmas.

\begin{lemma}
Let $y$ be the solution of equation \eqref{theequation} and $x\ge t$. Then
\begin{equation}\label{A:lemma:equality-yt=}
y(t)= y(x)\exp\Big(-\int_t^x q(s)ds\Big) + \int_t^x f(\xi) \exp\Big(-\int_t^\xi q(s)ds\Big)d\xi.
\end{equation}
\end{lemma}
\paragraph{Proof.}
From \eqref{theequation} it follows
\begin{equation*}
\frac{d}{d\xi} \biggl[ y(\xi) \exp\Big(-\int_0^\xi q(s)ds\Big)
\biggr] = -f(\xi) \exp\Big(-\int_0^\xi q(s)ds\Big),\quad \xi \in R.
\end{equation*}
To obtain \eqref{A:lemma:equality-yt=} one has to integrate the last equation from $t$ to $x$. \qed

\begin{lemma}
Let conditions I)-II) (see \S\ref{Preliminaries}) hold for some $p\in [1,\infty)$. Then
\begin{equation}\label{A:int_q=infty}
\int_{-\infty}^\infty q(\xi)d\xi = \infty.
\end{equation}
\end{lemma}
\paragraph{Proof.}
Assume the contrary, i.e., $\|q\|_1 < \infty$. Let $t\le x=0$ in \eqref{A:lemma:equality-yt=}. We obtain:
\begin{equation}\label{A:equality-yt-tlex}
y(t)=y(0)\exp\Big(-\int_t^0 q(s)ds\Big) + \int_t^0 f(\xi) \exp\Big(-\int_t^\xi q(s)ds\Big)d\xi.
\end{equation}
Let $f_0\in L_p(R) \cap L_1(R)$ and $f_0(\xi) >0,\ \xi\in R$ (for example,\ $f_0(\xi)= e^{-|\xi|}$). Then there exists $\lim\limits_{t\to -\infty } y(t)$. Indeed, the convergence of the first term in \eqref{A:equality-yt-tlex} is obvious. To check the convergence of the second term write it in the form:
\begin{equation}\label{A:equality-yt-tlex-secondterm}
\int_t^0 f_0(\xi) \exp\Big(-\int_t^\xi q(s)ds\Big)d\xi = \exp\Big(-\int_t^0 q(s)ds\Big)
\int_t^0 f_0(\xi) \exp\Big(\int^0_\xi q(s)ds\Big)d\xi.
\end{equation}
It is clear that the first factor on the right side of \eqref{A:equality-yt-tlex-secondterm} converges by our assumption, and the second one (by the same reason) converges as an improper integral:
\begin{equation*}
0\le \int_{-\infty}^0 f_0(\xi) \exp\Big( \int^0_\xi q(s)ds \Big) d\xi\ \le\ e^{\|q\|_1} \|f_0\|_1\ <\ \infty.
\end{equation*}
Hence there exists $A\eqdef \lim\limits_{t\to -\infty} y(t)$. Let us show that $A=0$. Indeed, if $A \ne 0$, then it is clear that
\begin{equation*}
\big|y(t)-A\big|\ \le\ \frac{|A|}{2}\quad \textrm{for all}\quad  t\ll -1,
\end{equation*}
and therefore,
\begin{equation*}
\big|y(t)\big| = \big|y(t) - A + A\big|\ \ge\ \Big|\big|A\big| -
\big|y(t)-A\big|\Big|\ \ge\ \frac{|A|}{2},\quad t\ll -1.
\end{equation*}
But from the last inequality it follows that $\|y\|_p = \infty$. Contradiction.\\
With $t$ tending to infinity in \eqref{A:equality-yt-tlex} we obtain
\begin{equation}\label{A:equality-yt-tlex-toinfinity}
0 = y(0)\exp\Big(-\int^0_{-\infty} q(s)ds \Big) + \int_{-\infty}^0 f_0(\xi)\exp\Big( -\int_{-\infty}^\xi q(s)ds \Big)d\xi.
\end{equation}
It is obvious that the second term in \eqref{A:equality-yt-tlex-toinfinity} is positive, since
\begin{equation*}
\int_{-\infty}^0 f_0(\xi)\exp\Big( -\int_{-\infty}^\xi
q(s)ds\Big)d\xi\ \ge\ \exp\Big( -\int_{-\infty}^0 q(s)ds \Big)
\int_{-\infty}^0 f_0(\xi)d\xi\ > \ 0;
\end{equation*}
and therefore, \eqref{A:equality-yt-tlex-toinfinity} implies $y(0) < 0$.\\

Now consider $x \ge t=0$ in \eqref{A:lemma:equality-yt=}:
\begin{equation}\label{A:equality-yt-xget}
y(0) = y(x)\exp\Big(-\int_0^x q(s)ds \Big) + \int^x_0 f_0(\xi)\exp\Big( -\int_0^\xi q(s)ds\Big)d\xi.
\end{equation}
Since we assumed that $\|q\|_1 < \infty$, from \eqref{A:equality-yt-xget} and by choosing $f_0$, it easily follows that there exists $\lim\limits_{x\to\infty} y(x)$ and it equals to $0$. Therefore, by passing to the limit in \eqref{A:equality-yt-xget} as $x \to\infty$, we obtain:
\begin{equation*}
y(0) =
\int_0^{\infty} f_0(\xi) \exp\Big( -\int_0^\xi q(s)ds \Big)d\xi \ge e^{-\|q\|_1} \int_0^{\infty} f_0(\xi) d\xi > 0.
\end{equation*}
Thus, $y(0) < 0 <y(0)$, contradiction. And finally the equality \eqref{A:int_q=infty} holds.\qed

We remind that in \S\ref{Preliminaries} the following function $d(x)$ was introduced:
\begin{equation}\label{A:definition_d}
d(x)=\inf_{d>0}\ \Big\{d : \int_{x-d}^{x+d}q(t)dt=2\Big\}.
\end{equation}
Because of \eqref{A:int_q=infty}, the function $d(x)$ is well-defined for all $x\in R$.

\begin{lemma}\label{A:lemma:sup_d}
Denote $d_0 \eqdef \sup\limits_{x\in R}\ d(x)$. Then
\begin{equation}\label{A:sup_d}
d_0 < \infty.
\end{equation}
\end{lemma}
\paragraph{Proof.}
Let $\psi(t) \in C_0^\infty(R)$,\ supp\ $\psi(t) = [-1,1],\ \psi(t) = 1$ for
$|t| \le \frac{1}{2}$, $0 \le \psi(t) \le 1$ for $t\in R$, i.e., $\psi(t)$ is a
cut-off function. Consider the function
\begin{equation}\label{A:z(t,x)}
z(t,x) = \psi\Big(\frac{t-x}{d(x)}\Big) \exp\Big(\int_x^t q(\xi)d\xi\Big),\quad t\in R,
\end{equation}
for fixed $x\in R$.\\
It is easy to see that $z(t,x)$ is the solution of the equation
\begin{equation}\label{A:equation_with_z}
-\frac{d}{dt}z(t,x) + q(t)z(t,x) = f(t,x),\quad t\in R,
\end{equation}
with the right side
\begin{equation}\label{A:rightside-equation_with_z}
f(t,x)= -\exp\Big(\int_x^t q(\xi)d\xi\Big)\
\Big[\psi\Big(\frac{t-x}{d(x)} \Big)\Big]'.
\end{equation}
To apply assertions I)-II) from \S\ref{Preliminaries} to \eqref{A:z(t,x)}, \eqref{A:equation_with_z} and \eqref{A:rightside-equation_with_z}, we need estimates of $\|f(\cdot,x)\|_p$ and
$\|z(\cdot,x)\|_p$. To this extent we use \eqref{A:definition_d} and properties of $\psi(t)$:
\begin{align}
\|f(\cdot,x)\|_p^p &= \int_{-\infty}^\infty |f(t,x)|^pdt =
\int_{x-d(x)}^{x+d(x)}\exp\Big(p\int_x^tq(\xi)d\xi\Big)
\bigg|\Big[\psi\Big(\frac{t-x}{d(x)}\Big)\Big]'\bigg|^p dt \nonumber \\
&\le\ c \exp \Big(p \int_{x-d(x)}^{x+d(x)} q(\xi)d\xi\Big)\
\frac{2d(x)}{d(x)^p} = \frac{2ce^{2p}}{d(x)^{p-1}}= \frac{c}{d(x)^{p-1}} \nonumber\\
&\Longrightarrow \qquad \|f(\cdot,x)\|_p\ \le\ cd(x)^{-1/p'},\quad p'=\frac{p}{p-1}.
\end{align}
Similarly, we obtain an upper estimate for $\|z(\cdot,x)\|_p$:
\begin{align}
\|z(\cdot,x)\|_p^p &= \int_{-\infty}^\infty |z(t,x)|^pdt =
\int_{x-d(x)}^{x+d(x)}\exp\Big(p\int_x^tq(\xi)d\xi\Big)
\Big|\psi\Big(\frac{t-x}{d(x)}\Big)\Big|^pdt\nonumber \\
&\le\ \exp \Big(p \int_{x-d(x)}^{x+d(x)} q(\xi)d\xi\Big)\ \ 2d(x) =
2e^{2p}\ d(x) = cd(x)\nonumber\\
&\Longrightarrow \qquad \|z(\cdot,x)\|_p\ \le\ cd(x)^{1/p}.
\end{align}
In addition, we can obtain a lower estimate for $\|z(\cdot,x)\|_p$:
\begin{align}\label{A:z:uppest-Lp}
\|z(\cdot,x)\|_p^p &=\int_{-\infty}^\infty |z(t,x)|^pdt\ \ge
\int_{x-d(x)/2}^{x+d(x)/2} |z(t,x)|^p dt \nonumber \\
&= \int_{x-d(x)/2}^{x+d(x)/2} \ \exp \Big(p \int_x^t q(\xi)d\xi\Big) \Big|\psi\Big(\frac{t-x}{d(x)}\Big)\Big|^p dt \nonumber \\
&\ge\ \exp \Big(-p \int_{x-d(x)}^{x+d(x)} q(\xi)d\xi\Big) \cdot d(x) = e^{-2p} d(x) = c^{-1} d(x)\nonumber \\
&\Longrightarrow \qquad \|z(\cdot,x)\|_p \ge  c^{-1} d(x)^{1/p}.
\end{align}
Since $d(x) \in (0, \infty)$ for every $x\in R$ (see \eqref{q}), the relations \eqref{A:equation_with_z}-\eqref{A:z:uppest-Lp} and I)-II) yield:
\begin{eqnarray*}
c^{-1} d(x)^{1/p}\ \le\ \|z(\cdot,x)\|_p\ \le\ c(p) \|f(\cdot,x)\|_p\ \le\ c d(x)^{-1/p'} \\
\Longrightarrow \ d(x)\le c,\ x\in R \quad \Longrightarrow \quad d_0=
\sup_{x\in R}\ d(x) \le c < \infty.
\end{eqnarray*}
\qed

From \eqref{A:sup_d} we obtain \eqref{A:inf_int_q} with $a=d_0$:
\begin{equation*}
q_o(d_0)= \inf_{x\in R} \int_{x-d_o}^{x+d_o} q(\xi)d\xi\
\ge\ \inf_{x\in R} \int_{x-d(x)}^{x+d(x)} q(\xi)d\xi = 2.
\end{equation*}

\paragraph{Proof of Theorem \ref{A:maintheorem}. \it{Sufficiency.}}
Let \eqref{A:inf_int_q} hold. Then
\begin{equation}\label{A:halfints_q=infty}
\int_{-\infty}^0 q(t)dt = \infty,\qquad \int^{\infty}_0 q(t)dt = \infty.
\end{equation}
The equalities \eqref{A:halfints_q=infty} can be checked in a similar way. Let us check the second one:
\begin{equation*}
\int_0^\infty q(t)dt = \sum_{k=0}^{\infty} \int_{2ka}^{2(k+1)a} q(t)dt\ \ge\ \sum_{k=0}^{\infty} q_0(a) = \infty.
\end{equation*}
From \eqref{A:inf_int_q} by Theorem \ref{bvp:theorem}, we obtain that in the space $L_p(R)$ for $p\in [1, \infty)$ the boundary value problem \eqref{bvp:equation}-\eqref{bvp:lim-y} is correctly solvable. Let $y$ be its solution. Then \eqref{theequation} has no other solutions in $L_p(R)$. Indeed, if $z$ is another solution of equation \eqref{theequation}, $z \ne y$ and $z \in L_p(R)$, then
\begin{equation*}
h(x) \eqdef y(x)-z(x) \in L_p(R).
\end{equation*}
In addition,
\begin{equation*}
\begin{split}
& -y'(x) + q(x)y(x) = f(x),\quad x\in R, \\
& -z'(x) + q(x)z(x) = f(x),\quad x\in R,
\end{split}
\end{equation*}
yield
\begin{equation}\label{A:anothersolution_h}
h(x) = c \exp \Big( \int_0^x q(\xi)d\xi\Big),\quad x\in R.
\end{equation}
But from \eqref{A:anothersolution_h} and \eqref{A:halfints_q=infty} it follows that $h\in L_p(R)$ only for $c=0$ and therefore $y(x) \equiv z(x),\ x\in R$, contradiction. By Theorem \ref{bvp:theorem} we obtain the statement of Theorem \ref{A:maintheorem}. In addition, since criterion \eqref{A:inf_int_q} does not depend on $p\in [1, \infty)$, only one of the statements A) or B) of Theorem \ref{A:maintheorem} is true. \qed

\subsection{Proof of main result for $p=\infty$.}

\paragraph{Proof of Theorem \ref{A:maintheorem}. \it{Necessity.}}

Let equation \eqref{theequation} be correctly solvable in $C(R)$. Then equality \eqref{A:int_q=infty} holds. Indeed, if on the contrary
\begin{equation*}
\int_{-\infty}^0 q(t)dt < \infty,\qquad \int^{\infty}_0 q(t)dt < \infty,
\end{equation*}
then let $y\in C(R)$ be a solution of equation \eqref{theequation} with $f\in C(R)$. Then the function
\begin{equation*}
z(x) = y(x) + \exp\Big(\int_0^x q(t)dt\Big),\quad x\in R,
\end{equation*}
is also a solution of equation \eqref{theequation} and $z\in C(R)$, which contradicts I)
(see \S\ref{Preliminaries}), and therefore \eqref{A:int_q=infty} holds. Let us check that $d_0<\infty$ (see \eqref{A:sup_d}). Below we will use \eqref{A:z(t,x)}, \eqref{A:equation_with_z} and \eqref{A:rightside-equation_with_z} from Lemma \ref{A:lemma:sup_d}. We need the following estimates of $\|f(\cdot,x)\|_{C(R)}$ and $\|z(\cdot,x)\|_{C(R)}$:
\begin{equation}\label{A:f:uppest-C}
\begin{split}
& \|f(\cdot,x)\|_{C(R)} = \sup_{t\in R} \exp\Big(\int_x^t q(\xi)d\xi\Big) \bigg|\Big[\psi\Big(\frac{t-x}{d(x)}\Big)\Big]'\bigg|\\
& \le\ c\exp\Big(\int_{x-d(x)}^{x+d(x)}q(\xi)d\xi\Big) \frac{1}{d(x)} = \frac{ce^2}{d(x)} = \frac{c}{d(x)};
\end{split}
\end{equation}
\begin{equation}\label{A:z:uppest-C}
\begin{split}
& \|z(\cdot,x)\|_{C(R)} = \sup_{t\in R}\ \psi\Big(\frac{t-x}{d(x)}\Big)\exp\Big(\int_x^t q(\xi)d\xi\Big) \\
&\le\ \exp\Big(\int_{x-d(x)}^{x+d(x)} q(\xi)d\xi\Big) = e^2 = c.
\end{split}
\end{equation}
In addition,
\begin{equation}\label{A:z:lowest-C}
\begin{split}
& \|z(\cdot,x)\|_{C(R)} = \sup_{t\in R}\ \psi\Big(\frac{t-x}{d(x)}\Big)
\exp\Big(\int_x^t q(\xi)d\xi\Big) \\
& \ge\ \exp\Big(- \int_{x-d(x)}^{x+d(x)} q(\xi)d\xi\Big)\ \sup_{t\in R}\
\psi\Big(\frac{t-x}{d(x)}\Big) = e^{-2} = c^{-1}.
\end{split}
\end{equation}
Since $d(x)\in (0, \infty)$ for $x\in R$ (see \eqref{q}), by \eqref{A:equation_with_z}-\eqref{A:rightside-equation_with_z}, \eqref{A:f:uppest-C}-\eqref{A:z:uppest-C}, I)-II) (see \S\ref{Preliminaries}) and \eqref{A:z:lowest-C}, we obtain:
\begin{equation*}
\begin{split}
& c^{-1}\ \le\ \|z(x)\|_{C(R)}\ \le\ c \|f(x)\|_{C(R)}\ \le\ \frac {c}{d(x)} \\
& \Rightarrow \quad d(x) \le c,\ x\in R \quad \Rightarrow \quad d_0 =
\sup_{x\in R}\ d(x) \le c < \infty.
\end{split}
\end{equation*}

But then \eqref{A:inf_int_q} holds because $q_0(d_0) \ge 2$ (see \S\ref{A:case_p<infty}).

\paragraph{Proof of Theorem \ref{A:maintheorem}. \it{Sufficiency.}}

Let $q_0(a) > 0$ for some $a\in (0, \infty)$. Then
\begin{equation}\label{A:J_0} 
J_o = \sup_{x\in R}\ J(x) < \infty,\quad J(x) \eqdef \int_x^{\infty}\exp\Big(-\int_x^t
q(\xi)d\xi\Big)dt,\quad x\in R.
\end{equation}

Indeed, \eqref{A:J_0} holds because of the following relations:
\begin{gather*}
J(x) = \int_x^{\infty} \exp \Big(- \int_x^t q(\xi)d\xi\Big)dt \\
= \int_x^{x+2a} \exp \Big(-\int_x^t q(\xi)d\xi \Big)dt\
+\ \sum_{k=1}^{\infty} \int_{x+2ka}^{x+2(k+1)a} \exp
\Big(-\int_x^t q(\xi)d\xi\Big)dt \\
\le\ 2a + \sum_{k=1}^{\infty}
2a \exp \Big(-\int_x^{x+2ka} q(\xi)d\xi\Big) \\
\le\ 2a + 2a \sum_{k=1}^{\infty} e^{-q_0(a)(k-1)}\ =\ 2a + 2a(1-e^{-q_0(a)})^{-1}\ <\ \infty.
\end{gather*}
Let $f\in C(R)$. Then it is obvious that the function
\begin{equation}\label{A:y(x)=}
y(x)=\int_x^\infty \exp\Big(-\int_x^t q(\xi)d\xi\Big) f(t)dt,\quad x\in R,
\end{equation}
is the solution of equation \eqref{theequation}, and
\begin{equation}\label{A:||y||_C}
\|y\|_{C(R)}\ \le\ \|f\|_{C(R)}\ \sup_{x\in R}\ J(x) = J_0\ \|f\|_{C(R)}.
\end{equation}

Under these conditions, the general solution $Y(x)$ of the equation \eqref{theequation} is of the form (see \eqref{A:y(x)=} and \eqref{A:anothersolution_h})
\begin{equation*}
Y(x) = y(x) + h(x),\quad x\in R;
\end{equation*}
and therefore, according to \eqref{A:||y||_C} and \eqref{A:int_q=infty}, we obtain $Y(x)\in C(R)$ only when $h(x)\equiv 0,\ x\in R$. Finally, statements I)-II) from \S\ref{Preliminaries} hold.
\qed

\subsection{Example}

\begin{example}
Let
\begin{equation}\label{A:example:q} 
q(x) = e^{x^2} + e^{x^2} \cos{e^{x^2}},\quad x\in R.
\end{equation}
Then equation \eqref{theequation} is correctly solvable in space $L_p(R),\ p\in [1,\infty]$.
\end{example}
\paragraph{Proof.}
Throughout below we denote
\begin{equation*}
q_1(x) = e^{x^2},\quad q_2(x) = e^{x^2} \cos{e^{x^2}},\quad x\in R.
\end{equation*}
Let $x\gg 1$. Consider the integral
\begin{equation*}
H(x) = \int_{x-1}^{x+1}\big(q_1(t) + q_2(t)\big)dt =
\int_{x-1}^{x+1}q_1(t)dt + \int_{x-1}^{x+1}q_2(t)dt, \quad x\in R.
\end{equation*}
Denote
\begin{equation*}
H_1(x) := \int_{x-1}^{x+1}q_1(t)dt,\ H_2(x) := \int_{x-1}^{x+1}q_2(t)dt.
\end{equation*}
Let us find estimates for $H_1(x)$ and $H_2(x)$, $x\gg 1$. We have:
\begin{equation*}
H_1(x) = \int_{x-1}^{x+1} e^{t^2}dt = \int_{x-1}^{x+1} \frac{de^{t^2}}{2t}
= \frac{e^{t^2}}{2t}\Big|^{x+1}_{x-1} + \frac{1}{2} \int_{x-1}^{x+1} \frac{e^{t^2}}{t^2}dt
\quad \Longrightarrow
\end{equation*}
\begin{equation*}
\int_{x-1}^{x+1} \Big(1- \frac{1}{2t^2}\Big) e^{t^2}dt
= \frac{e^{(x+1)^2}}{2(x+1)} - \frac{e^{(x-1)^2}}{2(x-1)} \quad \Longrightarrow
\end{equation*}
\begin{equation*}
\int_{x-1}^{x+1} e^{t^2}dt \ge  \int_{x-1}^{x+1} \Big(1- \frac{1}{2t^2}\Big)
e^{t^2}dt = \frac{e^{(x+1)^2}}{2(x+1)}\ \Big[1 - \frac{x+1}{x-1} e^{-4x}\Big]
\end{equation*}
\begin{equation*}
\ge \frac{e^{(x+1)^2}}{2(x+1)} \Big[ 1-2e^{-4x}\Big]  \ge \frac{e^{(x+1)^2}}{4(x+1)}.
\end{equation*}
Hence,
\begin{equation*}
H_1(x) \ge \frac{e^{(x+1)^2}}{4(x+1)},\quad x\gg 1.
\end{equation*}
Consider now $H_2(x),\ x\gg 1$:
\begin{equation*}
H_2(x) = \int_{x-1}^{x+1} e^{t^2}\cos{e^{t^2}}dt = \int_{x-1}^{x+1} \frac{d\sin{e^{t^2}}}{2t}
\end{equation*}
\begin{equation*}
= \frac{\sin{e^{t^2}}}{2t}\Big|_{x-1}^{x+1} + \frac{1}{2} \int_{x-1}^{x+1} \frac{\sin{e^{t^2}}}{t^2}dt \quad \Longrightarrow
\end{equation*}
\begin{equation*}
H_2(x) \le \frac{1}{2}\left( \frac{1}{x+1} + \frac{1}{x-1} \right) + \frac{1}{2} \left|\int_{x-1}^{x+1} \frac{\sin{e^{t^2}}}{t^2}dt \right|
\end{equation*}
\begin{equation*}
\le \frac{1}{2}\left( \frac{1}{x+1} + \frac{1}{x-1} \right) + \frac{1}{2} \int_{x-1}^{x+1} \frac{dt}{t^2} \le \frac{1}{x+1} + \frac{1}{x-1}
\end{equation*}
\begin{equation*}
= \frac{1}{x+1} \left[ 1 + \frac{x+1}{x-1} \right] \le \frac{4}{x+1}.
\end{equation*}
Hence,
\begin{equation*}
H_2(x) \le \frac{4}{x+1},\quad x\gg 1.
\end{equation*}
For $x\gg 1$ we find:
\begin{equation*}
H(x) = \int_{x-1}^{x+1} q(t)dt = H_1(x) + H_2(x)
\end{equation*}
\begin{equation*}
= H_1(x) \left[1+H_2(x)H_1^{-1}(x)\right] \ge H_1(x) \left[1-|H_2(x)|H_1^{-1}(x)\right]
\end{equation*}
\begin{equation*}
\ge H_1(x) \left[ 1 - \frac{4}{x+1}\cdot \frac{4(x+1)}{e^{(x+1)^2}} \right] =
H_1(x) \left[ 1 - 16 e^{-(x+1)^2}\right] \ge \frac{1}{2} H_1(x) \ge 1.
\end{equation*}
Thus, there exists $x_0 \gg 1$ such that
\begin{equation*}
H(x) \ge 1, \quad x\ge x_0.
\end{equation*}
Note here that the function $q$ is even and we get finally:
\begin{equation}\label{A:example:Hge1-|x|gex0} 
H(x) \ge 1, \quad |x|\ge x_0.
\end{equation}
Now consider the function $H(x)$ for $x\in [-x_0,x_0]$. Let $a=2x_0+2$.  Then

\begin{eqnarray*}
\left.
\begin{array}{c}
x-a \le x_0 - a = -x_0 -2 \\
x+a \ge -x_0 + a = x_0 + 2
\end{array}
\right.
\quad \Longrightarrow \\
\end{eqnarray*}
\begin{equation}\label{A:example:union_included}
[-x_0-2,-x_0] \cup[x_0,x_0+2]\subseteq [x-a,x+a].
\end{equation}
\\
From \eqref{A:example:Hge1-|x|gex0} and \eqref{A:example:union_included} it follows that
\begin{eqnarray}\label{A:example:int_q-ge2}
\int_{x-a}^{x+a}q(t)dt &\ge& \int_{-x_0-2}^{-x_0}q(t)dt + \int_{x_0}^{x_0+2}q(t)dt \nonumber \\
&=& \int_{(-x_0-1)-1}^{(-x_0-1)+1}q(t)dt+\int_{(x_0+1)-1}^{(x_0+1)+1}q(t)dt \nonumber \\
&=& H(-x_0-1)+H(x_0+1)\ \ge\ 2.
\end{eqnarray}
At the same time, for $|x|\le x_0$ from \eqref{A:example:Hge1-|x|gex0} and by choosing $a=2x_0+2$ it easily follows that
\begin{equation}\label{A:example:int_q-ge2-|x|gex0}
\int_{x-a}^{x+a}q(t)dt \ge 2, \quad x\in [-x_0, x_0].
\end{equation}
Thus, by \eqref{A:example:Hge1-|x|gex0}, \eqref{A:example:int_q-ge2} and \eqref{A:example:int_q-ge2-|x|gex0} we conclude that there exists $a>0$ such that
\begin{equation*}
\int_{x-a}^{x+a}q(t)dt \ge 2,\quad x\in R.
\end{equation*}
Then the equation \eqref{theequation} is correctly solvable in space $L_p(R),\ p\in [1,\infty]$ by Theorem \ref{A:maintheorem}.


\section{Solution of Problem B}

\subsection{Sharp by order estimates for one class of definite integrals}\label{B1:proof} 
\setcounter{equation}{0}

\paragraph{}
In this section we consider the integral
\begin{equation*}
J(x) = \int_{x}^{\infty} e^{-\int\limits_{x}^{t}q(\xi)d\xi}dt,\quad x\in R,
\end{equation*}
with non-negative locally integrated function $q$. Our goal is to find sharp by order two-sided estimates for the function $J$. In particular, below we prove the following statement:

\begin{theorem}\label{B1:maintheorem}
Let condition \eqref{q} for $q$ hold. Suppose that there exist an absolutely continuous function $q_1(x)>0$ and a function $q_2\in L_1^{loc}(R)$ such that
\begin{equation*}
q(x) = q_1(x) + q_2(x),\quad x\in R.
\end{equation*}
Let there exist a continuous and positive for all $x\in R$ function $s(x)$ such that the following conditions hold:
\begin{enumerate}
\item[a)]
\begin{equation*}
s(x) \to \infty \quad \textrm{as} \quad |x|\to \infty,
\end{equation*}
\item[b)]
\begin{equation*}
\frac{1}{s(x)} \ge \frac{|q'_1(x)|}{q_1^2(x)} \quad \textrm{for all} \quad |x| \gg 1,
\end{equation*}
\item[c)]
\begin{equation*}
\lim_{|x|\to \infty} \frac{s(x)}{xq_1(x)}=0,
\end{equation*}
\item[d)]
for some $\nu\in[1,\infty)$ the inequalities hold:
\begin{equation*}
\frac{1}{\nu}\le\frac{s(t)}{s(x)}\le\nu ,\quad t\in \Delta(x),\ |x|\gg 1,
\end{equation*}
where
\begin{equation*}
\Delta(x)=\big[\Delta^{-}(x),\ \Delta^{+}(x)\big]=\Big[x-\frac{s(x)}{q_1(x)},\ x+\frac{s(x)}{q_1(x)}\Big],\quad x\in R,
\end{equation*}
\item[e)]
\begin{equation*}
\varkappa_0 = \sup_{x\in R} |\tilde{\varkappa}(x)| < \infty,
\end{equation*}
where
\begin{equation*}
\tilde{\varkappa}(x) = \sup_{t\in \Delta(x)} |\varkappa(t)|,
\end{equation*}
\begin{equation*}
\varkappa(t) = q_1(t) \int_{x}^{t} \frac{q_2(\xi)}{q_1(\xi)}d\xi,\quad t\in \Delta(x),\ x\in R.
\end{equation*}
\end{enumerate}
Then the following estimates hold:
\begin{equation}\label{B1:mainresult}
\frac{c^{-1}}{q_1(x)} \le J(x)\le \frac{c}{q_1(x)},\quad x\in R.
\end{equation}
\end{theorem}
\paragraph{Proof.}
Let $x\in R$. Denote by $\{\tilde{\Omega}_{n}\}_{n=1}^{\infty}$ segments of the real axis such that
\begin{eqnarray}\label{B1:DefinitionBOmega_n} 
\left\{
\begin{array}{ll}
\displaystyle
\tilde{\Omega}_n=\tilde{\Omega}(x_n) =\big[\tilde{\Omega}_n^{-},\ \tilde{\Omega}_n^{+}\big]= 
\Big[x_n,\ x_n + \frac{s(x_n)}{q_1(x_n)}\Big],\ n\ge 1,\\
\displaystyle
x_{n+1} = x_n + \frac{s(x_n)}{q_1(x_n)},\ n\ge 1,\\
x_1 = x.
\end{array}
\right.
\end{eqnarray}
From \eqref{B1:DefinitionBOmega_n} it follows:
\begin{equation}\label{B1:UnionBOmega_n}
\bigcup_{n=1}^{\infty} \tilde{\Omega}_n = [x, \infty)
\end{equation} 
Indeed, if \eqref{B1:UnionBOmega_n} does not hold, then there exists  $z\in (x, \infty)$ such that
\begin{equation}\label{B1:ContrUnionOmega_n}
\tilde{\Omega}_n^+ = \tilde{\Omega}^+(x_n) < z, \quad n\ge 1.
\end{equation}
Then, obviously, $x_n < z$ for $n\ge 1$. Since a sequence $\{x_n\}_{n=1}^{\infty}$ is monotonically ascending (by construction) and bounded, it has a limit $x_0 \in (x, z]$. Note here that \eqref{B1:DefinitionBOmega_n} and \eqref{B1:ContrUnionOmega_n} yield:
\begin{eqnarray*}
\infty &>& z-x\ \ge\ \sum_{n=1}^{\infty} (\tilde{\Omega}_n^+ - \tilde{\Omega}_n^-) = \sum_{n=1}^{\infty} \frac{s(x_n)}{q_1(x_n)} \quad \Longrightarrow \\
&&0= \lim_{n\to \infty} \frac{s(x_n)}{q_1(x_n)} = \frac{s(x_0)}{q_1(x_0)},
\end{eqnarray*}
contradiction, because functions $s$ and $q_1$ are positive. Hence,  \eqref{B1:UnionBOmega_n} holds.

\begin{remark}
{\rm The proof of \eqref{B1:UnionBOmega_n} is obtained by Otelbayev' method (see \cite{MO}).}
\end{remark}

The following equalities are based on \eqref{B1:DefinitionBOmega_n}-\eqref{B1:UnionBOmega_n}:
\begin{equation}\label{B1:J=sum} 
J(x) = \int_{x}^{\infty} e^{-\int\limits_{x}^{t}q(\xi)d\xi}dt =
\sum_{n=1}^{\infty}\int_{\tilde{\Omega}_n}e^{-\int\limits_{x_1}^t q(\xi)d\xi}dt,\quad x\in R.
\end{equation}

Below to obtain estimates for terms in \eqref{B1:J=sum} we need various auxiliary statements. We emphasize here that all conditions of Theorem \ref{B1:maintheorem} are assumed to hold.

\begin{lemma}\label{B1:lemma:int_q-split} 
Let $t\in [x, \Delta^+(x)],\ x\in R$. Then
\begin{equation}\label{B1:int_q-split}
\int_{x}^{t}q(\xi)d\xi = \int_{x}^{t}q_1(\xi)d\xi + \varkappa(t) + \int_{x}^{t}q_1(\xi)\mu(\xi)d\xi,
\end{equation}
where
\begin{equation*}
\mu(t) = -\ \frac{q'_1(t)}{q^2_1(t)}\ \varkappa(t).
\end{equation*}
\end{lemma}
\paragraph{Proof.}
To check the equality \eqref{B1:int_q-split} we use an integration by parts:
\begin{align*}
\int_{x}^{t} q(\xi)d\xi &= \int_{x}^{t} q_1(\xi)d\xi + \int_{x}^{t} q_2(\xi)d\xi \\
&= \int_{x}^{t} q_1(\xi)d\xi + \int_{x}^{t} q_1(\xi) \frac{q_2(\xi)}{q_1(\xi)} d\xi \\
&= \int_{x}^{t} q_1(\xi)d\xi + \int_{x}^{t} q_1(\xi)\ d\bigg[\int_{x}^{\xi} \frac{q_2(\tau)}{q_1(\tau)} d\tau \bigg] \\
&= \int_{x}^{t} q_1(\xi)d\xi +\varkappa(t) - \int_{x}^{t} q'_1(\xi) \left( \int_{x}^{\xi} \frac{q_2(\tau)}{q_1(\tau)} d\tau \right) d\xi \\
&= \int_{x}^{t} q_1(\xi)d\xi +\varkappa(t) + \int_{x}^{t} q_1(\xi) \left[ -\ \frac{q'_1(\xi)}{q^2_1(\xi)}\ q_1(\xi) \left( \int_{x}^{\xi} \frac{q_2(\tau)}{q_1(\tau)} d\tau \right) \right] d\xi \\
&= \int_{x}^{t} q_1(\xi)d\xi + \varkappa(t) + \int_{x}^{t} q_1(\xi)\mu(\xi)d\xi.
\end{align*}
\qed

\begin{lemma}{\rm\cite{Maj}}\label{B1:lemma:q1t/q1x}
Let $|x| \gg 1$. Then
\begin{equation}\label{B1:q1t/q1x} 
e^{-\ \frac{\nu}{s(x)} \int\limits_{x}^{t} q_1(\xi)d\xi} \le \frac{q_1(t)}{q_1(x)} \le e^{\frac{\nu}{s(x)} \int\limits_{x}^{t} q_1(\xi)d\xi},\quad t\in [x, \Delta^+(x)].
\end{equation}
\end{lemma}
\paragraph{Proof.}
From conditions of Theorem \ref{B1:maintheorem} for $\xi\in [x, \Delta^+(x)]$ we get relations:
\begin{eqnarray*}
-\ \nu\ \frac{q_1(\xi)}{s(x)}\ \le\ -\ \frac{q_1(\xi)}{s(\xi)}\ \le\ \frac{q'_1(\xi)}{q_1(\xi)}\ \le\ \frac{q_1(\xi)}{s(\xi)}\ \le\ \nu\ \frac{q_1(\xi)}{s(x)} \quad \Longrightarrow \\
-\ \frac{\nu}{s(x)} \int_{x}^{t} q_1(\xi)d\xi\ \le\ \ln\frac{q_1(t)}{q_1(x)}\ \le\ \frac{\nu}{s(x)} \int_{x}^{t} q_1(\xi)d\xi.
\end{eqnarray*}
\qed

\begin{lemma}{\rm\cite{Maj}}\label{B1:lemma:sln/nu}
For $|x|\gg 1$ the following inequality holds:
\begin{equation}\label{B1:sln/nu} 
\frac{s(x)}{\nu} \ln(1+\nu)\ \le\ \int_{x}^{\Delta^+(x)} q_1(\xi)d\xi.
\end{equation}
\end{lemma}
\paragraph{Proof.}
From \eqref{B1:q1t/q1x} we obtain:
\begin{eqnarray*}
\begin{array}{c}
\displaystyle
\frac{\nu q_1(x)}{s(x)} \le\ \frac{\nu q_1(t)}{s(x)}\
e^{\frac{\nu}{s(x)}\int\limits_{x}^{t}q_1(\xi)d\xi},\quad t\in [x, \Delta^+(x)]\ \Longrightarrow\\
\displaystyle
\nu =\ \frac{\nu q_1(x)}{s(x)}\ \Big[\Delta^+(x)-x\Big] \le \int_{x}^{\Delta^+(x)} \frac{\nu q_1(t)}{s(x)}\ e^{\frac{\nu}{s(x)}\int\limits_{x}^{t}q_1(\xi)d\xi} dt \\
= e^{\frac{\nu}{s(x)}\int\limits_{x}^{t}q_1(\xi)d\xi}
\Big|_{x}^{\Delta^+(x)} = \ e^{\frac{\nu}{s(x)}\int\limits_{x}^{\Delta^+(x)}q_1(\xi)d\xi} - 1\ \Longrightarrow\\
1+\nu \le e^{\frac{\nu}{s(x)}\int\limits_{x}^{\Delta^+(x)}q_1(\xi)d\xi}\ \Longrightarrow \quad \eqref{B1:sln/nu}
\end{array}
\end{eqnarray*}
\qed

\begin{lemma}\label{B1:two-inequalities}
For $|x| \gg 1$ the following inequalities hold:
\begin{equation}\label{B1:int_Omega-q-ge} 
\int_{x}^{\Delta^+(x)} q(\xi)d\xi \ge \frac{1}{4} \int_{x}^{\Delta^+(x)} q_1(\xi)d\xi,
\end{equation}
\begin{equation}\label{B1:int_xt-inOmega-q-ge} 
\int_x^t q(\xi)d\xi \ge \frac{1}{4} \int_x^t q_1(\xi)d\xi -\varkappa_0,\quad t\in [x, \Delta^+(x)].
\end{equation}
\end{lemma}
\paragraph{Proof.}
To obtain \eqref{B1:int_Omega-q-ge} we use Lemma \ref{B1:lemma:int_q-split}. We note that for all $|x| \gg 1$ and $\xi\in [x, \Delta^+(x)]$ the following inequalities hold:
\begin{equation*}
|\mu(\xi)|=\frac{|q'_1(\xi)|}{q^2_1(\xi)}|\varkappa(\xi)| \le \frac{\varkappa_0}{s(\xi)}\le \frac{1}{2}.
\end{equation*}
According to \eqref{B1:int_q-split} we have:
\begin{equation*}
\int_x^{\Delta^+(x)}q(\xi)d\xi \ge \int_x^{\Delta^+(x)}q_1(\xi)d\xi - \varkappa_0 - \frac{1}{2} \int_x^{\Delta^+(x)}q_1(\xi)d\xi = \frac{1}{2} \int_x^{\Delta^+(x)}q_1(\xi)d\xi -\varkappa_0,
\end{equation*}
and by \eqref{B1:sln/nu} we finally obtain
\begin{eqnarray*}
\int_{x}^{\Delta^+(x)}q(\xi)d\xi \ge \frac{1}{2} \int_{x}^{\Delta^+(x)}q_1(\xi)d\xi \bigg[1 - 2\varkappa_0 \bigg(\int_{x}^{\Delta^+(x)}q_1(\xi)d\xi\bigg)^{-1}\bigg] \\
\ge \frac{1}{2} \int_{x}^{\Delta^+(x)}q_1(\xi)d\xi \bigg[1- \frac{2\varkappa_0\ \nu}{\ln(1+\nu)} \frac{1}{s(x)}\bigg] \ge \frac{1}{4} \int_{x}^{\Delta^+(x)} q_1(\xi)d\xi.
\end{eqnarray*}
The inequality \eqref{B1:int_xt-inOmega-q-ge} is obtained in the same way. \qed

\begin{lemma} 
Let $x\gg1$ and $t\ge x$. Then
\begin{equation}\label{B1:int_xt-q-ge}
\int_{x}^{t} q(\xi)d\xi \ge \frac{1}{4} \int_x^t q_1(\xi)d\xi -\varkappa_0.
\end{equation}
\end{lemma}
\paragraph{Proof.}
Let $\{\tilde{\Omega}_n\}_{n=1}^{\infty}$ be segments of \eqref{B1:DefinitionBOmega_n}. It is clear that there exists $n \ge 1$ such that $t \in \tilde{\Omega}_n$ (see \eqref{B1:UnionBOmega_n}). If $n=1$ then \eqref{B1:int_xt-q-ge} is equivalent to \eqref{B1:int_xt-inOmega-q-ge}. If $n>1$ then \eqref{B1:DefinitionBOmega_n}, \eqref{B1:int_Omega-q-ge}, and \eqref{B1:int_xt-inOmega-q-ge} yield
\begin{eqnarray*}
\int_x^t q(\xi)d\xi &=& \sum_{k=1}^{n-1} \int_{\tilde{\Omega}_k} q(\xi)d\xi + \int_{x_n}^{t} q(\xi)d\xi \\
&\ge& \sum_{k=1}^{n-1} \frac{1}{4} \int_{\tilde{\Omega}_k} q_1(\xi)d\xi+ \frac{1}{4}\int_{x_n}^{t}  q_1(\xi)d\xi - \varkappa_0 \\
&=& \frac{1}{4} \int_{x}^{t} q_1(\xi)d\xi -\varkappa_0.
\end{eqnarray*}
\qed

\begin{corollary}
Let $x\gg 1$. Then
\begin{equation}\label{B1:JleJ1} 
J(x) \le c\ J_1(x),\quad J_1(x) =
\int_{x}^{\infty} e^{-\frac{1}{4} \int\limits_{x}^{t}q_1(\xi)d\xi}dt.
\end{equation}
\end{corollary}
\paragraph{Proof.}
According to \eqref{B1:int_xt-q-ge} we have:
\begin{equation*}
J(x) = \int_{x}^{\infty} e^{- \int\limits_{x}^{t}q(\xi)d\xi}dt \le e^{\varkappa_0} \int_{x}^{\infty} e^{-\frac{1}{4} \int\limits_{x}^{t}q_1(\xi)d\xi}dt.
\end{equation*}
\qed

\begin{lemma}
Let $x\gg 1$. Then
\begin{equation}\label{B1:J1-upest} 
J_1(x) \le \frac{8}{q_1(x)}.
\end{equation}
\end{lemma}
\paragraph{Proof.}
Let $a\ge 0$. Denote
\begin{equation*}
J_1(x,a) = \int_{x}^{x+a} e^{-\frac{1}{4} \int\limits_{x}^{t}q_1(\xi)d\xi}dt.
\end{equation*}
The following relations are based on integration by parts:
\begin{eqnarray*}
J_1(x,a) &=& \int_{x}^{x+a} \Big(-\frac{4}{q_1(t)}\Big) \Big(-\frac{1}{4}q_1(t)\Big)\ e^{-\frac{1}{4} \int\limits_{x}^{t}q_1(\xi)d\xi}dt\\
&=& - \frac{4}{q_1(t)}\ e^{-\frac{1}{4} \int\limits_{x}^{t}q_1(\xi)d\xi} \Big|^{x+a}_{x} - \int_{x}^{x+a} \frac{4q'_1(t)}{q^2_1(t)}\ e^{-\frac{1}{4} \int\limits_{x}^{t}q_1(\xi)d\xi} dt\\
&\le& \frac{4}{q_1(x)} - \frac{4}{q_1(x+a)}\ e^{-\frac{1}{4} \int\limits_{x}^{x+a}q_1(\xi)d\xi} + \int_{x}^{x+a} \frac{4}{s(t)}\ e^{-\frac{1}{4} \int\limits_{x}^{t}q_1(\xi)d\xi} dt\\
&\le& \frac{4}{q_1(x)} + \frac{1}{2} \int_{x}^{x+a} e^{-\frac{1}{4} \int\limits_{x}^{t}q_1(\xi)d\xi} dt = \frac{4}{q_1(x)} + \frac{1}{2} J_1(x,a)\\
&&\qquad \Longrightarrow \quad J_1(x,a) \le \frac{8}{q_1(x)}, \quad x\gg 1.
\end{eqnarray*}
By the above inequality we obtain \eqref{B1:J1-upest}:
\begin{equation*}
J_1(x) = \lim_{a\to \infty} J_1(x,a) \le \frac{8}{q_1(x)}.
\end{equation*}
\qed \\
Now, for $x\gg 1$  the upper estimate in \eqref{B1:mainresult} follows from \eqref{B1:JleJ1} and \eqref{B1:J1-upest}. Hence, there exists $x_0\gg 1$ such that for $x\ge x_0$ we have:
\begin{equation}\label{B1:J-xgex0upest}
J(x) \le \frac{c(x_0)}{q_1(x)}
\end{equation}
Let now $x\le -x_0$ (see \eqref{B1:J-xgex0upest}). Then
\begin{equation*}
J(x) = \int_{x}^{\infty} e^{-\int\limits_{x}^{t}q(\xi)d\xi} dt =
\int_{x}^{-x_0} e^{-\int\limits_{x}^{t}q(\xi)d\xi} dt + \int_{-x_0}^{\infty} e^{-\int\limits_{x}^{t}q(\xi)d\xi} dt.
\end{equation*}
Denote
\begin{equation*}
W(x,-x_0) := \int_{x}^{-x_0} e^{-\int\limits_{x}^{t}q(\xi)d\xi} dt, \quad x\le -x_0.
\end{equation*}
Then
\begin{equation}\label{B1:J=wxx} 
J(x) = W(x,-x_0) \left[ 1+\frac{1}{W(x,-x_0)}  \int_{-x_0}^{\infty}
e^{-\int\limits_{x}^{t}q(\xi)d\xi} dt \right].
\end{equation}
Notice here that
\begin{equation*}
\frac{1}{W(x,-x_0)} \int_{-x_0}^{\infty}e^{-\int\limits_{x}^{t}q(\xi)d\xi}dt =
\frac{e^{-\int\limits_{x}^{-x_0}q(\xi)d\xi}}{W(x,-x_0)} \int_{-x_0}^{\infty} e^{-\int\limits_{-x_0}^{t}q(\xi)d\xi}dt.
\end{equation*}
It is easy to see that by \eqref{B1:J-xgex0upest} the integral
\begin{equation*}
\int_{-x_0}^{\infty} e^{-\int\limits_{-x_0}^{t}q(\xi)d\xi}dt\ :=\ A(x_0),
\end{equation*}
converges since
\begin{equation*}
\int_{-x_0}^{\infty} e^{-\int\limits_{-x_0}^{t}q(\xi)d\xi}dt =
\int_{-x_0}^{x_0} e^{-\int\limits_{-x_0}^{t}q(\xi)d\xi}dt \ + \ e^{-\int\limits_{-x_0}^{x_0}q(\xi)d\xi} J(x_0).
\end{equation*}
From the other side,
\begin{eqnarray*}
e^{\int\limits_{x}^{-x_0}q(\xi)d\xi}\cdot W(x,-x_0) \ = \
e^{\int\limits_{x}^{-x_0}q(\xi)d\xi} \int_{x}^{-x_0} e^{-\int\limits_{x}^{t}q(\xi)d\xi}dt \\
= \int_{x}^{-x_0} e^{\int\limits_{t}^{-x_0}q(\xi)d\xi}dt\ \ge\ (|x| - x_0) \to \infty
\quad \textrm{as} \quad x \to -\infty.
\end{eqnarray*}
Therefore, there exists $x_1\gg x_0$ such that for $x\le -x_1$ by \eqref{B1:J=wxx} it follows that
\begin{equation}\label{B1:Jle3w}
J(x) \le 3W(x,-x_0),\quad x\le -x_1.
\end{equation}
Let us now find an estimate for $W(x,-x_0)$, $x\le -x_1$.
Let $\{\tilde\omega_{n}\}_{n=-\infty}^{-1}$ be segments, constructed by the following rule:
\begin{eqnarray}\label{B1:DefinitionSOmega_n} 
\left\{
\begin{array}{ll}
\displaystyle
\tilde\omega_n = \tilde\omega(x_n) = \big[\tilde\omega_n^{-},\ \tilde\omega_n^{+}\big] = 
\Big[x_n - \frac{s(x_n)}{q_1(x_n)},\ x_n\Big],\ n\le -1,\\
\displaystyle
x_{n-1} = x_n - \frac{s(x_n)}{q_1(x_n)},\ n\le -1,\\
x_{-1}= -x_0.
\end{array}
\right.
\end{eqnarray}
From \eqref{B1:DefinitionSOmega_n} the equality follows:
\begin{equation*}
\bigcup_{n=-\infty}^{-1} \tilde\omega_n = (-\infty, -x_0],
\end{equation*}
which is checked in the same way as \eqref{B1:UnionBOmega_n}.\\
Below we need the analogies of lemmas we obtained before, which can be proved by the very same methods.

\begin{lemma}
Let $t\in [\Delta^-(x),x]$. Then (see \eqref{B1:int_q-split})
\begin{equation*}
\int^{x}_{t}q(\xi)d\xi = \int^{x}_{t}q_1(\xi)d\xi - \varkappa(t) - \int^{x}_{t}q_1(\xi)\mu(\xi)d\xi.
\end{equation*}
\end{lemma}

\begin{lemma}{\rm\cite{Maj}} 
Let $|x| \gg 1$. Then
\begin{equation*}
e^{-\ \frac{\nu}{s(x)} \int\limits_{t}^{x} q_1(\xi)d\xi} \ \le \ \frac{q_1(t)}{q_1(x)} \ \le \ e^{\frac{\nu}{s(x)} \int\limits_{t}^{x} q_1(\xi)d\xi},\quad t\in [\Delta^-(x),x].
\end{equation*}
\end{lemma}

\begin{lemma}{\rm\cite{Maj}}
For $|x|\gg 1$ the following inequality holds:
\begin{equation*}
\frac{s(x)}{\nu} \ln(1+\nu)\ \le\ \int_{\Delta^-(x)}^{x} q_1(\xi)d\xi
\end{equation*}
\end{lemma}

\begin{lemma}
For all $x \le -x_0$ the following inequalities hold:
\begin{equation*}
\int_{\Delta^-(x)}^{x} q(\xi)d\xi \ \ge \ \frac{1}{4} \int_{\Delta^-(x)}^{x} q_1(\xi)d\xi,
\end{equation*}
\begin{equation*}
\int_{\Delta^-(x)}^t q(\xi)d\xi \ \ge \ \frac{1}{4} \int_{\Delta^-(x)}^t q_1(\xi)d\xi -\varkappa_0, \quad t\in [\Delta^-(x),x].
\end{equation*}
\end{lemma}

\begin{lemma}
Let $x\le t\le -x_0$. Then
\begin{equation*}
\int_{x}^{t} q(\xi)d\xi \ \ge \ \frac{1}{4} \int_x^t q_1(\xi)d\xi -\varkappa_0.
\end{equation*}
\end{lemma}

\begin{corollary}
Let $x\le -x_0$. Then
\begin{equation}\label{B1:wlew1}
W(x,-x_0) \le c\ W_1(x,-x_0),\quad W_1(x,-x_0) = \int_{x}^{-x_0}
e^{-\frac{1}{4} \int\limits_{x}^{t}q_1(\xi)d\xi}dt.
\end{equation}
\end{corollary}

\begin{lemma}
Let $x\le -x_0$. Then
\begin{equation}\label{B1:w1-upest}
W_1(x,-x_0) \le \frac{8}{q_1(x)}.
\end{equation}
\end{lemma}
By \eqref{B1:Jle3w}, \eqref{B1:wlew1} and \eqref{B1:w1-upest} we obtain that for $x\le -x_1$ the following inequalities hold:
\begin{equation}\label{B1:J-xlex1upest}
J(x) \ \le \ 3W(x,-x_0) \ \le \ 3c\,W_1(x,-x_0) \ \le \ \frac{24c}{q_1(x)}
\end{equation}
From \eqref{B1:J-xgex0upest} and \eqref{B1:J-xlex1upest} it follows that
\begin{equation}\label{B1:J1-notmidupest}
J(x) \le \frac{c}{q_1(x)}, \quad x\notin (-x_1, x_0).
\end{equation}
Now consider the function $\varphi(x) = q_1(x) J(x)$ for $x\in [-x_1, x_0]$. Since the integral $J(x)$ converges (by our proof above), the function $\varphi(x)$ is continuous on $[-x_1, x_0]$  and therefore has on this interval a finite maximum, i.e.
\begin{equation}\label{B1:J1-midupest} 
q_1(x)J(x)\le c < \infty, \quad x\in [-x_1, x_0].
\end{equation}
Then the upper estimate of \eqref{B1:mainresult} follows directly from \eqref{B1:J1-notmidupest} and \eqref{B1:J1-midupest}.
We obtain now the lower estimate in \eqref{B1:mainresult}. The following base inequality holds for any $x\in R$.
\begin{equation*}
J(x) = \int_{x}^{\infty} e^{-\int\limits_{x}^{t}q(\xi)d\xi}dt \ge \int_{x}^{\Delta^+(x)} e^{-\int\limits_{x}^{t}q(\xi)d\xi}dt,\quad x\in R.
\end{equation*}
Let $|x|\gg 1,\ t\in [x,\Delta^+(x)]$. The following relations follow from Lemma \ref{B1:lemma:int_q-split}:
\begin{eqnarray*}
\int_{x}^{t}q(\xi)d\xi &=& \int_{x}^{t}q_1(\xi)d\xi
+ \varkappa(t) + \int_{x}^{t}q_1(\xi)\mu(\xi)d\xi \\
&\le& \frac{3}{2} \int_{x}^{t}q_1(\xi)d\xi + \varkappa_0 \quad \Longrightarrow \\
J(x) &\ge& \int_{x}^{\Delta^+(x)} e^{-\int\limits_{x}^{t}q(\xi)d\xi}dt \ge \int_{x}^{\Delta^+(x)} e^{-\frac{3}{2}\int\limits_{x}^{t}q_1(\xi)d\xi}e^{-\varkappa_0}dt \\
&=& e^{-\varkappa_0} \int_{x}^{\Delta^+(x)} e^{-\frac{3}{2}\int\limits_{x}^{t}q_1(\xi)d\xi}dt.
\end{eqnarray*}
Denote
\begin{equation*}
\tilde{J_1}(x)= \int_{x}^{\Delta^+(x)}e^{-\frac{3}{2}
\int\limits_{x}^{t}q_1(\xi)d\xi}dt,\quad x\in R.
\end{equation*}
Let $|x|\gg 1$. Then
\begin{gather*}
\tilde{J_1}(x) = -\frac{2}{3} \int_{x}^{\Delta^+(x)} \frac{1}{q_1(t)} \left[-\frac{3}{2} q_1(t) e^{-\frac{3}{2}\int\limits_{x}^{t}q_1(\xi)d\xi}\right] dt\\
= \frac{2}{3q_1(x)} - \frac{2}{3}\frac{1}{q_1(\Delta^+(x))}\ e^{-\frac{3}{2}\int\limits_{x}^{\Delta^{+}(x)}q_1(\xi)d\xi} - \frac{2}{3} \int_{x}^{\Delta^+(x)} \frac{q'_1(t)}{q^2_1(t)}\ e^{-\frac{3}{2}\int\limits_{x}^{t}q_1(\xi)d\xi}dt \\
\ge \frac{2}{3q_1(x)} - \frac{2}{3q_1(x)} \frac{q_1(x)}{q_1(\Delta^+(x))}\ e^{-\frac{3}{2}\int\limits_{x}^{\Delta^{+}(x)}q_1(\xi)d\xi} - \frac{2}{3}  \int_{x}^{\Delta^+(x)} \frac{1}{s(t)}\ e^{-\frac{3}{2}\int\limits_{x}^{t}q_1(\xi)d\xi}dt \\
\ge \frac{2}{3q_1(x)} - \frac{2}{3q_1(x)}\ e^{-(\frac{3}{2} - \frac{\nu}{s(x)}) \int\limits_x^{\Delta^+(x)} q_1(\xi)d\xi} -\ \frac{1}{2}\tilde{J_1}(x)\quad \Longrightarrow \\ \\
\frac{3}{2} \tilde{J_1}(x) \ge \frac{2}{3q_1(x)} - \frac{2}{3q_1(x)}\
e^{-\int\limits_{x}^{\Delta^+(x)} q_1(\xi)d\xi} 
\ge \frac{2}{3q_1(x)} - \frac{2}{3q_1(x)}\ e^{-\frac{s(x)}{\nu}\ln{(1+\nu)}} \\
= \frac{2}{3q_1(x)} \Big[1 - \frac{1}{(1+\nu)^{\frac{s(x)}{\nu}}}\Big] \ge \frac{1}{3q_1(x)} \quad \Longrightarrow
\end{gather*}
\begin{equation}\label{B1:J-notmidlowest} 
J(x) \ge c^{-1}\tilde{J_1}(x) \ge \frac{c^{-1}}{q_1(x)}, \quad |x|\gg 1.
\end{equation}
Thus, the inequality \eqref{B1:J-notmidlowest} holds for $|x|\ge \tau\gg 1$. \\
Let $\varphi(x) = q_1(x) J(x),\ x\in [-\tau, \tau]$. The functions $q_1(x)$ and $J(x)$ are both positive for $x\in R$  and continuous. Therefore $\varphi(x)$ has on $[-\tau, \tau]$ a positive minimum $\varphi_0$. Hence,
\begin{equation*}
q_1(x)J(x) = \varphi(x) \ge \varphi_0,\quad x\in [-\tau, \tau] \Longrightarrow
\end{equation*}
\begin{equation}\label{B1:J-midlowest} 
J(x) \ge \frac{c^{-1}}{q_1(x)},\quad x\in [-\tau, \tau]
\end{equation}
The lower estimate in \eqref{B1:mainresult} follows from \eqref{B1:J-notmidlowest} and \eqref{B1:J-midlowest}. Theorem \ref{B1:maintheorem} is proved. \qed

\subsection{The asymptotic formula for calculation of one class of definite integrals}\label{B2:proof} 

In this section we continue to study the integral
\begin{equation*}
J(x) = \int_{x}^{\infty} e^{-\int\limits_{x}^{t}q(\xi)d\xi}dt,\quad x\in R,
\end{equation*}
with non-negative locally integrated function $q$. Our goal here is to make inequalities \eqref{B1:mainresult} more accurate. The main result of this section is the following Theorem \ref{B2:maintheorem}:
\begin{theorem}\label{B2:maintheorem}
Let condition \eqref{q} for $q$ hold. Suppose that there exist an absolutely continuous function $q_1(x)>0$ and a function $q_2\in L_1^{loc}(R)$ such that
\begin{equation*}
q(x) = q_1(x) + q_2(x),\quad x\in R.
\end{equation*}
Let there exist a continuous and positive for all $x$ function $s(x)$ such that the following conditions hold:
\begin{enumerate}
\item[a)]
\begin{equation*}
s(x) \to \infty \quad \textrm{as} \quad |x|\to \infty,
\end{equation*}
\item[b)]
\begin{equation*}
\frac{1}{s(x)} \ge \frac{|q'_1(x)|}{q_1^2(x)} \quad \textrm{for all} \quad |x| \gg 1,
\end{equation*}
\item[c)]
\begin{equation*}
\lim_{|x|\to \infty} \frac{s(x)}{xq_1(x)}=0,
\end{equation*}
\item[d)]
for some $\nu\in[1,\infty)$ the inequalities hold:
\begin{equation*}
\frac{1}{\nu}\le\frac{s(t)}{s(x)}\le\nu,\quad t\in \Delta(x),\ |x|\gg 1,
\end{equation*}
where
\begin{equation*}
\Delta(x)=\big[\Delta^{-}(x),\ \Delta^{+}(x)\big]=\Big[x-\frac{s(x)}{q_1(x)},\ x+\frac{s(x)}{q_1(x)}\Big],\quad x\in R,
\end{equation*}
\item[e)]
\begin{equation*}
\tilde{\varkappa}(x) \to 0 \quad \textrm{as} \ |x| \to \infty,
\end{equation*}
where
\begin{equation*}
\tilde{\varkappa}(x) = \sup_{t\in \Delta(x)} |\varkappa(t)|, \quad x\in R,
\end{equation*}
\begin{equation*}
\varkappa(t) = q_1(t) \int_{x}^{t} \frac{q_2(\xi)}{q_1(\xi)}d\xi,\quad t\in \Delta(x).
\end{equation*}
\end{enumerate}
Then
\begin{equation}\label{B2:mainresult}
J(x) = \int_{x}^{\infty} e^{-\int\limits_{x}^{t}q(\xi)d\xi}dt =
\frac{1+\varepsilon(x)}{q_1(x)}, \quad \lim_{|x|\to \infty} \varepsilon(x) = 0.
\end{equation}
\end{theorem}

Below for the proof of Theorem \ref{B2:maintheorem} we need various auxiliary statements. We emphasize that all conditions of Theorem \ref{B2:maintheorem} are assumed to hold.

\paragraph{Proof.}
Denote the integral $\tilde{J}(x)$ and the function $\tilde{\mu}(x)$:
\begin{gather*}
\tilde{J}(x) = \int_{x}^{\Delta^+(x)} e^{-\int\limits_{x}^{t}q(\xi)d\xi}dt,\quad x\in R,\\
\tilde{\mu}(x)=\sup_{t\in [x,\Delta^+(x)]}|\mu(t)|,\quad \mu(t)=-\frac{q'_1(t)}{q^2_1(t)}\varkappa(t),\quad t\in [x,\Delta^+(x)].
\end{gather*}

\begin{lemma}\label{B2:lemma:ge-int_xt-q-le} 
Let $|x|\gg 1$, $t\in [x,\Delta^+(x)]$. Then
\begin{equation}\label{B2:ge-int_xt-q-le} 
(1-\tilde{\mu}(x)) \int_{x}^{t}q_1(\xi)d\xi - \tilde{\varkappa}(x) \le \int_{x}^{t}q(\xi)d\xi \le
(1+\tilde{\mu}(x)) \int_{x}^{t}q_1(\xi)d\xi + \tilde{\varkappa}(x).
\end{equation}
\end{lemma}
\paragraph{Proof.}
Let $t\in [x,\Delta^+(x)]$. Then (see \eqref{B1:int_q-split})
\begin{equation*}
\int_{x}^{t}q(\xi)d\xi = \int_{x}^{t}q_1(\xi)d\xi + \varkappa(t) + \int_{x}^{t}q_1(\xi)\mu(\xi)d\xi.
\end{equation*}
From this equality we obtain \eqref{B2:ge-int_xt-q-le} directly:
\begin{eqnarray*}
\int_{x}^{t}q(\xi)d\xi &\ge&
\int_{x}^{t}q_1(\xi)d\xi - \tilde{\varkappa}(x) - \tilde{\mu}(x) \int_{x}^{t}q_1(\xi)d\xi \\
&=& (1-\tilde{\mu}(x)) \int_{x}^{t}q_1(\xi)d\xi - \tilde{\varkappa}(x),\\
\int_{x}^{t}q(\xi)d\xi &\le&
\int_{x}^{t}q_1(\xi)d\xi + \tilde{\varkappa}(x) + \tilde{\mu}(x) \int_{x}^{t}q_1(\xi)d\xi \\
&=& (1+\tilde{\mu}(x)) \int_{x}^{t}q_1(\xi)d\xi + \tilde{\varkappa}(x).
\end{eqnarray*}
\qed

\begin{corollary}
Let $|x|\gg 1$. Then
\begin{equation}\label{B2:ge-Jtilde-le}
\begin{array}{c}
\displaystyle
\tilde{J}(x) \ge
(1-c\tilde{\varkappa}(x))\int_{x}^{\Delta^+(x)}
e^{-(1+\tilde{\mu}(x))\int\limits_{x}^{t}q_1(\xi)d\xi} dt,\\
\tilde{J}(x) \le
\displaystyle
(1+c\tilde{\varkappa}(x))\int_{x}^{\Delta^+(x)}
e^{-(1-\tilde{\mu}(x))\int\limits_{x}^{t}q_1(\xi)d\xi} dt.
\end{array}
\end{equation}
\end{corollary}
\paragraph{Proof.}
Since $\Delta^+(x) \to \pm \infty$ as $x\to \pm \infty$ (see condition {\it c)}), we have \\
$\tilde{\varkappa}(x) \to 0$,\ \ $\tilde{\mu}(x) \to 0$ as $|x| \to \infty$.
Then
\begin{equation*}
e^{\tilde{\varkappa}(x)} \lessgtr (1\pm c\tilde{\varkappa}(x)), \quad |x| \to \infty.
\end{equation*}
Therefore, for $|x|\gg 1$ and $t\in [x,\Delta^+(x)]$ we obtain:
\begin{eqnarray*}
&e^{-\int\limits_{x}^{t}q(\xi)d\xi}& \le
e^{-(1-\tilde{\mu}(x)) \int\limits_{x}^{t}q_1(\xi)d\xi} e^{\tilde{\varkappa}(x)} \le
(1+c\tilde{\varkappa}(x)) e^{-(1-\tilde{\mu}(x)) \int\limits_{x}^{t}q_1(\xi)d\xi},\\
&e^{-\int\limits_{x}^{t}q(\xi)d\xi}& \ge
e^{-(1+\tilde{\mu}(x)) \int\limits_{x}^{t}q_1(\xi)d\xi} e^{-\tilde{\varkappa}(x)} \ge
(1-c\tilde{\varkappa}(x)) e^{-(1+\tilde{\mu}(x)) \int\limits_{x}^{t}q_1(\xi)d\xi}.
\end{eqnarray*}
The above inequalities easily imply relations in \eqref{B2:ge-Jtilde-le}, for example
\begin{equation*}
\tilde{J}(x) = \int_{x}^{\Delta^+(x)} e^{-\int\limits_{x}^{t}q(\xi)d\xi}dt \le
(1+c\tilde{\varkappa}(x)) \int_{x}^{\Delta^+(x)}
e^{-(1-\tilde{\mu}(x)) \int\limits_{x}^{t}q_1(\xi)d\xi}dt.
\end{equation*}
The second inequality in \eqref{B2:ge-Jtilde-le} is obtained in the same way. \qed

\begin{lemma}
Fix $\alpha > 0$. Then for $|x|\gg 1$ the equality holds:
\begin{equation}\label{B2:Jtildealpha-asymptotics}
\tilde{J}_{\alpha}(x) = \frac{1}{\alpha q_1(x)} \left(1+ O\Big(\frac{1}{s(x)}\Big) \right),
\end{equation}
where
\begin{equation*}
\tilde{J}_{\alpha}(x) = \int_{x}^{\Delta^+(x)}e^{-\alpha \int\limits_{x}^{t}q_1(\xi)d\xi}dt,
\quad x\in R.
\end{equation*}
\end{lemma}
\paragraph{Proof.}
Below we use Lemmas \ref{B1:lemma:q1t/q1x} and \ref{B1:lemma:sln/nu}:
\begin{eqnarray*}
\tilde{J}_{\alpha}(x) &=& \int_{x}^{\Delta^+(x)} e^{-\alpha \int\limits_{x}^{t} q_1(\xi)d\xi}dt \\ &=& \frac{1}{q_1(x)} \int_{x}^{\Delta^+(x)} \frac{q_1(x)}{q_1(t)}
q_1(t) e^{-\alpha \int\limits_{x}^{t} q_1(\xi)d\xi} dt \\
&\le& \frac{1}{q_1(x)} \int_{x}^{\Delta^+(x)} q_1(t) e^{-(\alpha - \frac{\nu}{s(x)}) \int\limits_{x}^{t} q_1(\xi)d\xi} dt \\
&=& \frac{1}{\frac{\nu}{s(x)}-\alpha} \frac{1}{q_1(x)} \left[ e^{(\frac{\nu}{s(x)}-\alpha) \int\limits_{x}^{t} q_1(\xi)d\xi}\Big\vert_{x}^{\Delta^+(x)} \right] \\
&=& \frac{1}{\alpha-\frac{\nu}{s(x)}} \frac{1}{q_1(x)} \left[ 1 - e^{(\frac{\nu}{s(x)} -\alpha) \int\limits_{x}^{\Delta^+(x)} q_1(\xi)d\xi} \right] \\
&\le& \frac{1}{\alpha q_1(x)} \left[ 1+ O\Big(\frac{1}{s(x)}\Big) \right];
\end{eqnarray*}
likewise,
\begin{eqnarray*}
\tilde{J}_{\alpha}(x) &=& \int_{x}^{\Delta^+(x)} e^{-\alpha \int\limits_{x}^{t} q_1(\xi)d\xi} dt \\
&=& \frac{1}{q_1(x)} \int_{x}^{\Delta^+(x)} \frac{q_1(x)}{q_1(t)} q_1(t) e^{-\alpha \int\limits_{x}^{t} q_1(\xi)d\xi} dt \\
&\ge& \frac{1}{q_1(x)} \int_{x}^{\Delta^+(x)} q_1(t) e^{-(\alpha + \frac{\nu}{s(x)}) \int\limits_{x}^{t} q_1(\xi)d\xi} dt \\
&=& \frac{-1}{\frac{\nu}{s(x)}+\alpha} \frac{1}{q_1(x)} \left[ e^{-(\frac{\nu}{s(x)}+\alpha) \int\limits_{x}^{t} q_1(\xi)d\xi}\Big\vert_{x}^{\Delta^+(x)} \right] \\
&=& \frac{1}{\alpha+\frac{\nu}{s(x)}} \frac{1}{q_1(x)} \left[ 1- e^{-(\alpha +\frac{\nu}{s(x)}) \int\limits_{x}^{\Delta^+(x)} q_1(\xi)d\xi} \right] \\
&\ge& \frac{1}{\alpha+\frac{\nu}{s(x)}} \frac{1}{q_1(x)} \left[1- \frac{1}{(1+\nu)^{1+\frac{\alpha s(x)}{\nu}}} \right].
\end{eqnarray*}
Note here that
\begin{equation*}
(1+\nu)^{1+\frac{\alpha s(x)}{\nu}} \ge (1+\nu)^{\frac{\alpha s(x)}{\nu}} \ge 1+ \nu \frac{\alpha s(x)}{\nu} = 1+ \alpha s(x) \ge \alpha s(x),
\end{equation*}
and therefore,
\begin{equation*}
\tilde{J}_{\alpha}(x) \ge \frac{1}{\alpha + \frac{\nu}{s(x)}} \frac{1}{q_1(x)}
\left[1+ O\Big(\frac{1}{s(x)}\Big)\right]
= \frac{1}{\alpha q_1(x)} \left[ 1+ O\Big(\frac{1}{s(x)}\Big) \right].
\end{equation*}
The above inequalities yield \eqref{B2:Jtildealpha-asymptotics}. \qed

\begin{corollary}\label{B2:cor:Jtilde=}
Let $|x|\gg 1$. Then
\begin{equation*}
\tilde{J}(x) = \int_{x}^{\Delta^+(x)} e^{-\int\limits_{x}^{t}q(\xi)d\xi}dt = \frac{1+\delta(x)}{q_1(x)}, 
\end{equation*}
where
\begin{equation*}
|\delta(x)| \le c \left( \frac{1}{s(x)} + \tilde{\varkappa}(x) \right).
\end{equation*}
\end{corollary}
\paragraph{Proof.}
The corollary follows from \eqref{B2:ge-Jtilde-le}, \eqref{B2:Jtildealpha-asymptotics} and obvious estimates
\begin{equation*}
\tilde{\mu}(x) = \sup_{t\in [x,\Delta^+(x)]} |\mu(t)| = \sup_{t\in [x,\Delta^+(x)]} \frac{|q'_1(t)|}{q^2_1(t)} |\varkappa(t)| \le \sup_{t\in [x,\Delta^+(x)]} \frac{\tilde\varkappa(x)}{s(t)} \le \frac{\nu \tilde\varkappa(x)}{s(x)}
\end{equation*}
\qed
\\
Now we can finish a proof of \eqref{B2:mainresult}. Below we use Theorem \ref{B1:maintheorem}, Corollary~\ref{B2:cor:Jtilde=}, and Lemmas \ref{B1:lemma:q1t/q1x} and \ref{B1:lemma:sln/nu}:
\begin{eqnarray*}
J(x) &=& \int_{x}^{\infty} e^{-\int\limits_{x}^{t}q(\xi)d\xi}dt =
\int_{x}^{\Delta^+(x)} e^{-\int\limits_{x}^{t}q(\xi)d\xi}dt +
\int^{\infty}_{\Delta^+(x)} e^{-\int\limits_{x}^{t}q(\xi)d\xi}dt\\
&=& \tilde{J}(x) + e^{-\int\limits_{x}^{\Delta^+(x)}q(\xi)d\xi} J(\Delta^+(x))\\
&=& \tilde{J}(x) \Big[ 1 + e^{-\int\limits_{x}^{\Delta^+(x)}q(\xi)d\xi} J(\Delta^+(x)) \tilde{J}^{-1}(x) \Big] := \tilde{J}(x) [1 + \eta(x)].
\end{eqnarray*}
But by Lemma \ref{B2:lemma:ge-int_xt-q-le} we have:
\begin{eqnarray*}
|\eta(x)| &\le& e^{-\int\limits_{\Delta(x)}q(\xi)d\xi} \frac{c}{q_1(\Delta^{+}(x))}\ q_1(x) \\
&\le& c\ e^{\tilde{\varkappa}(x) - (1-\tilde{\mu}(x)) \int\limits_{\Delta(x)} q_1(\xi)d\xi} \frac{q_1(x)}{q_1(\Delta^+(x))} \le
c\ e^{-\frac{1}{2} \int\limits_{\Delta(x)}q_1(\xi)d\xi} \frac{q_1(x)}{q_1(\Delta^+(x))} \\
&\le& c\ e^{-\frac{1}{2} \int\limits_{\Delta(x)}q_1(\xi)d\xi} \cdot
e^{\frac{\nu}{s(x)}\int\limits_{\Delta(x)}q_1(\xi)d\xi} =
c^{-(\frac{1}{2} - \frac{\nu}{s(x)}) \int\limits_{\Delta(x)}q_1(\xi)d\xi} \\
&\le& c\ e^{-\frac{1}{4} \int\limits_{\Delta(x)}q_1(\xi)d\xi} \le c\ e^{-\frac{1}{4} \frac{s(x)}{\nu}\ln{(1+\nu)}}\le \frac{c}{s(x)},
\end{eqnarray*}
and, finally,
\begin{equation*}
J(x) = \tilde{J}(x) \left( 1 + O\Big(\frac{1}{s(x)}\Big) \right) =
\frac{1 + O(\frac{1}{s(x)}) + \tilde{\varkappa}(x)}{q_1(x)}.
\end{equation*}
\qed

\subsection{The Theorem on asymptotic majorant}\label{B3:proof} 

In this section we prove Theorem \ref{B:maintheorem} (see \S\ref{B}). For convenience we remind below its statement.
\\ \\
{\bf Theorem \ref{B:maintheorem}.}
{\it
Let condition \eqref{q} for $q$ hold. Suppose that there exist an absolutely continuous function $q_1(x)>0$ and a function $q_2\in L_1^{loc}(R)$ such that
\begin{equation*}
q(x) = q_1(x) + q_2(x),\quad x\in R.
\end{equation*}
Let there exist a continuous and positive for all $x$ function $s(x)$ such that the following conditions hold:\\
\begin{enumerate}
\item[a)]
\begin{equation*}
s(x) \to \infty \quad \textrm{as} \quad |x|\to \infty,
\end{equation*}
\item[b)]
\begin{equation}\label{B:condition_b}
\frac{1}{s(x)} \ge \frac{|q'_1(x)|}{q_1^2(x)} \quad \textrm{for all} \quad |x| \gg 1,
\end{equation}
\item[c)]
\begin{equation}\label{B:condition_c}
\lim_{|x|\to \infty} \frac{s(x)}{xq_1(x)}=0,
\end{equation}
\item[d)]
for some $\nu\in[1,\infty)$ the inequalities hold:
\begin{equation}\label{B:condition_d}
\frac{1}{\nu}\le\frac{s(t)}{s(x)}\le\nu,\quad t\in \Delta(x),\ |x|\gg 1,
\end{equation}
where
\begin{equation*}
\Delta(x)=\Big[\Delta^{-}(x),\ \Delta^{+}(x)\Big]=\Big[x-\frac{s(x)}{q_1(x)},\ x+\frac{s(x)}{q_1(x)}\Big],\quad x\in R,
\end{equation*}
\item[e)]
\begin{equation}\label{B:condition_e}
\tilde{\varkappa}(x) \to 0 \ \textrm{ as } \ |x| \to \infty,
\end{equation}
where
\begin{equation*}
\tilde{\varkappa}(x) = \sup_{t\in \Delta(x)} |\varkappa(t)|, \quad x\in R,
\end{equation*}
\begin{equation*}
\varkappa(t) = q_1(t) \int_{x}^{t} \frac{q_2(\xi)}{q_1(\xi)}d\xi,\quad t\in \Delta(x).
\end{equation*}
\end{enumerate}
Then for $p\in [1, \infty]$ an asymptotic majorant $\varkappa_p(x)$ for solutions $y\in D_p$
\begin{equation*}
D_p = \big\{\ y: -y'(x) + q(x)y(x) = f(x), \ y\in L_p(R),\ \|f\|_p=1\ \big\}
\end{equation*}
is of the form:
\begin{eqnarray}
\varkappa_p(x)= 
\left\{
\begin{array}{cc}
1, & \quad p=1,\ x\in R, \\
\displaystyle \frac{1}{(p')^{1/p'}} \frac{1}{q_1(x)^{1/p'}}, & \quad
p\in (1,\infty),\ p'=\frac{p}{p-1},\ |x|\gg 1,
\vspace{2mm}\\
\displaystyle \frac{1}{q_1(x)}, & \quad p=\infty, \ |x|\gg 1.
\end{array}
\right.
\end{eqnarray}
}

\paragraph{Proof of Theorem \ref{B:maintheorem}.}
Fix $x\in R$. Then a value of solution $y\in D_p$ in the point $x$ is of the form (see \eqref{y=int-G(x,t)fdt}-\eqref{G(x,t)=}):
\begin{equation}\label{B3:y=int-G(x,t)fdt}
y(x) \eqdef (Gf)(x) = \int_{-\infty}^{\infty} G(x,t)f(t)dt,\quad x\in R,\\
\end{equation}
where
\begin{equation}\label{B3:G(x,t)=}
G(x,t) =
\left\{
\begin{array}{l}
0, \quad  t<x,\ x\in R,\\
\displaystyle
\exp \Big( -\int_{x}^{t}q(\xi)d\xi \Big),\quad t\ge x,\ x\in R.
\end{array}
\right.
\end{equation}
It is clear that the formula \eqref{B3:y=int-G(x,t)fdt}-\eqref{B3:G(x,t)=} represents a linear functional $G$, defined on $L_p(R)$. By general properties of linear functionals, defined on~$L_p(R)$ (see \cite[ch.V, \S2.2-2.3, ch.VI. \S2)]{KA}, we obtain:
\begin{eqnarray}\label{B3:||G||=formula}
\|G(x)\| &=& \sup_{y\in D_p} |y(x)| =
\left\{
\begin{array}{ll}
\esssup\limits_{t\in R} G(x,t), & p=1, \\
\bigg( \displaystyle \int_{-\infty}^{\infty}G(x,t)^{p'}dt \bigg)^{1/p'},
& p\in (1,\infty),\ p'=\frac{p}{p-1},
\vspace{1mm}\\
\displaystyle \int_{-\infty}^{\infty}G(x,t)dt, & p=\infty.
\end{array}
\right.
\nonumber \\
&&= \
\left\{
\begin{array}{ll}
\esssup\limits_{t\ge x}\ e^{-\int\limits_{x}^{t}q(\xi)d\xi}, & p=1, \\
\bigg( \displaystyle \int_{x}^{\infty} e^{-p'\int\limits_{x}^{t}q(\xi)d\xi}dt \bigg)^{1/p'},
& p\in (1,\infty),\ p'=\frac{p}{p-1}, \\
\displaystyle \int_{x}^{\infty}e^{-\int\limits_{x}^{t}q(\xi)d\xi}dt, & p=\infty.
\end{array}
\right.
\end{eqnarray}
It is clear that for $|x|\gg 1$ by Theorem \ref{B2:maintheorem} we obtain from \eqref{B3:||G||=formula} that
\begin{eqnarray}\label{B3:||G||=result}
\|G(x)\| &=& \Big(1+\varepsilon(x)\Big) \cdot
\left\{
\begin{array}{cc}
1, & p=1\\
\displaystyle \frac{1}{(p')^{1/p'}} \frac{1}{q_1(x)^{1/p'}}, & p\in (1,\infty)
\vspace{2mm}\\
\displaystyle \frac{1}{q_1(x)}, & p=\infty
\end{array}
\right. , \\ \nonumber
\textrm{ where } &&
\lim\limits_{|x|\to \infty} \varepsilon(x) = 0.
\end{eqnarray}
Now a statement of Theorem \ref{B:maintheorem} follows from definition of asymptotic majorant $\varkappa_p(x)$, relations \eqref{B3:||G||=formula}-\eqref{B3:||G||=result}, and definition of supremum.

\qed

\subsection{Example}

\begin{example}
Let
\begin{equation}\label{B:example:q} 
q(x) = e^{x^2} + e^{x^2} \cos{e^{x^2}},\quad x\in R.
\end{equation}
Let $\varkappa_p(x)$ be an asymptotic majorant of solutions for equation \eqref{theequation}. Then
\begin{eqnarray}\label{B:example:Teta}
\varkappa_p(x) = 
\left\{
\begin{array}{ll}
\quad \quad 1, & \quad p=1,\ x\in R, \\
\displaystyle \frac{1}{(p')^{1/p'}} \frac{1}{e^{x^2/p'}}, & \quad
p\in (1,\infty),\ |x|\gg 1, \\
\displaystyle \quad \ \ \frac{1}{e^{x^2}}, & \quad p=\infty,\ |x|\gg 1.
\end{array}
\right.
\end{eqnarray}
\end{example}
\paragraph{Proof.}
Let us show that if
\begin{equation}\label{B:example:q1q2}
q_1(x) = e^{x^2},\quad q_2(x) = e^{x^2}\cos{e^{x^2}},\quad x\in R,
\end{equation}
and
\begin{equation}\label{B:example:s(x)}
s(x) = \frac{e^{x^2}}{8\sqrt{1+x^2}},\quad x\in R,
\end{equation}
then all conditions of Theorem \ref{B:maintheorem} hold and therefore \eqref{B:example:Teta} follows from it directly. We have to check only relations \eqref{B:condition_b}, \eqref{B:condition_c}, \eqref{B:condition_d}, and \eqref{B:condition_e}. By choosing \eqref{B:example:q1q2}-\eqref{B:example:s(x)} all other conditions of Theorem \ref{B:maintheorem} obviously hold. We have:
\begin{equation*}
\lim_{|x|\to \infty} = \frac{s(x)}{xq_1(x)}
= \lim_{|x|\to \infty} \frac{1}{8x\sqrt{1+x^2}} = 0\quad
\Longrightarrow \quad \eqref{B:condition_c}
\end{equation*}
Now, since $1+a^2\ge 2a$,
\begin{equation*}
\frac{1}{s(x)} = \frac{8\sqrt{1+x^2}}{e^{x^2}}\ge \frac{2|x|}{e^{x^2}} =
\frac{|q'_1(x)|}{q^2_1(x)}\quad \Longrightarrow \quad \eqref{B:condition_b}
\end{equation*}
Let us obtain \eqref{B:condition_d}. We first calculate $\Delta(x)$:
\begin{equation*}
\Delta(x)=\Big[ x-\frac{s(x)}{q_1(x)},\ x+\frac{s(x)}{q_1(x)}\Big]=\Big[x-\frac{1}{8\sqrt{1+x^2}},\ x+\frac{1}{8\sqrt{1+x^2}} \Big]
\end{equation*}
Let $t\in [x,\Delta^+(x)]$. Then according to Lagrange' formula we obtain:
\begin{equation*}
s(t) = s(x) + s'(\xi)(t-x),\quad x\in R,
\end{equation*}
where $\xi$ is between $t$ and $x$.\\
Let
\begin{equation*}
M(x) = \max_{t\in [x,\Delta^+(x)]}\frac{s(t)}{s(x)},\quad m(x)=\min_{t\in [x,\Delta^+(x)]}\frac{s(t)}{s(x)}.
\end{equation*}
Next,
\begin{equation*}
\frac{|s'(x)|}{s(x)} = \Big|2x - \frac{x}{1+x^2}\Big| = |x|\Big|2-\frac{1}{1+x^2}\Big|
\le 2|x|,\quad x\in R.
\end{equation*}
Therefore,
\begin{equation*}
M(x) = \max_{t\in [x,\Delta^+(x)]} \Big| 1 + \frac{s'(\xi)}{s(x)}(t-x) \Big| \le
1+ \max_{t\in [x,\Delta^+(x)]} \frac{|s'(\xi)|}{s(\xi)}\cdot \frac{s(\xi)}{s(x)}\frac{1}{8\sqrt{1+x^2}}
\end{equation*}
\begin{equation*}
\le 1 + \frac{M(x)}{8\sqrt{1+x^2}} \max_{\xi\in [x,\Delta^+(x)]} 2|\xi| \le
1 +  \frac{M(x)}{4\sqrt{1+x^2}} \Big( |x| + \frac{1}{8\sqrt{1+x^2}} \Big)
\end{equation*}
\begin{equation*}
= 1 + \frac{M(x)}{4} \Big( \frac{|x|}{\sqrt{1+x^2}} + \frac{1}{8(1+x^2)}\Big) \le 1 + \frac{M(x)}{4} \Big(1+\frac{1}{8}\Big) = 1 + \frac{9}{32}M(x) \quad \Longrightarrow
\end{equation*}
\begin{equation}\label{B:example:Mle}
M(x) \le \frac{32}{23},\quad x\in R.
\end{equation}
Now we use \eqref{B:example:Mle} to obtain an estimate for $m(x)$:
\begin{equation*}
m(x) = \min_{t\in [x,\Delta^+(x)]} \Big| 1 + \frac{s'(\xi)}{s(x)}(t-x) \Big| \ge
\min_{t\in [x,\Delta^+(x)]} \bigg[ 1 - \frac{s'(\xi)}{s(\xi)}\cdot \frac{s(\xi)}{s(x)}|t-x| \bigg]
\end{equation*}
\begin{equation*}
\ge 1 - \frac{M(x)}{8\sqrt{1+x^2}} \max_{\xi\in [x,\Delta^+(x)]} 2|\xi| \ge 1 - \frac{8}{23}
\left[ \frac{|x|}{\sqrt{1+x^2}}+\frac{1}{8(1+x^2)} \right]
\end{equation*}
\begin{equation*}
\ge 1 - \frac{8}{23}\Big(1+\frac{1}{8}\Big) = \frac{14}{23} > \frac{1}{2}\quad \Longrightarrow
\end{equation*}
\begin{equation*}
\frac{1}{2}\le m(x) \le \frac{s(t)}{s(x)} \le M(x) \le 2,\quad t\in [x,\Delta^+(x)],\quad x\in R.
\end{equation*}
For $t\in [\Delta^-(x),x]$ the proof of \eqref{B:condition_d} goes along the same lines. Finally,
\begin{equation*}
\frac{1}{2}\le m(x) \le \frac{s(t)}{s(x)} \le M(x) \le 2,\quad t\in \Delta(x),\quad x\in R.
\end{equation*}
It remained to check \eqref{B:condition_e}. Let $t\in [x,\Delta^+(x)]$. We have:
\begin{equation}\label{B:example:kappa}
\varkappa(t) = q_1(t) \int_{x}^{t} \frac{q_2(\xi)}{q_1(\xi)}d\xi
= e^{t^2} \int_{x}^{t} \cos{e^{\xi^2}}d\xi,\quad t\in [x,\Delta^+(x)].
\end{equation}
From \eqref{B:example:kappa} it follows that
\begin{equation}\label{B:example:kappa'}
\varkappa'(t) = 2t e^{t^2} \int_{x}^{t} \cos{e^{\xi^2}}d\xi + e^{t^2}\cos{e^{t^2}},\quad t\in [x,\Delta^+(x)].
\end{equation}
Hence, according to \eqref{B:example:kappa}-\eqref{B:example:kappa'}, the function $\varkappa(t)$ is a solution of Cauchy' problem:
\begin{eqnarray}\label{B:example:cauchy_problem}
\left\{
\begin{array}{c}
\varkappa'(t) = 2t\varkappa(t) + e^{t^2}\cos{e^{t^2}},\quad t \in [x,\Delta^+(x)], \\
\varkappa(x) = 0.
\end{array}
\right.
\end{eqnarray}
Relations in \eqref{B:example:cauchy_problem} yield
\begin{equation}\label{B:example:kappa=}
\varkappa(t) = \int_{x}^{t}2\xi\varkappa(\xi)d\xi + \int_{x}^{t}e^{\xi^2}\cos{e^{\xi^2}}d\xi, \quad t\in [x,\Delta^+(x)].
\end{equation}
Let $\Phi(x) = \max\limits_{t\in [x,\Delta^+(x)]} |\varkappa(t)|$.
Since $\Phi(x)=|\varkappa(x_0)|$ for some $x_0\in [x,\Delta^+(x)]$, by \eqref{B:example:kappa=} we obtain:
\begin{equation*}
\Phi(x) \le \max_{t\in [x,\Delta^+(x)]} \left| \int_{x}^{t}2\xi \varkappa(\xi)d\xi \right| + \max_{t\in [x,\Delta^+(x)]} \left| \int_{x}^{t}e^{\xi^2}\cos{e^{\xi^2}d\xi} \right|
\end{equation*}
\begin{equation*}
\le 2\Phi(x) \max_{t\in [x,\Delta^+(x)]} \left| \int_{x}^{t}\xi d\xi \right|
+ \max_{t\in [x,\Delta^+(x)]} \left| \int_{x}^{t}e^{\xi^2}\cos{e^{\xi^2}d\xi} \right|
\end{equation*}
\begin{equation*}
= \Phi(x) \max_{t\in [x,\Delta^+(x)]}|t^2-x^2|+ \max_{t\in [x,\Delta^+(x)]}
\left| \int_{x}^{t}e^{\xi^2}\cos{e^{\xi^2}d\xi} \right|
\end{equation*}
But for $t\in [x,\Delta^+(x)]$ we have:
\begin{equation*}
\big|t^2-x^2\big|=\big|t-x\big|\big|t+x\big|\le \frac{1}{8\sqrt{1+x^2}}\left[ 2|x| + \frac{1}{8\sqrt{1+x^2}} \right]
\end{equation*}
\begin{equation*}
\le \frac{|x|}{4\sqrt{1+x^2}} + \frac{1}{64(1+x^2)}\le \frac{1}{2}\quad \Longrightarrow
\end{equation*}
\begin{equation*}
\Phi(x) \le \frac{1}{2}\Phi(x) + \max_{t\in [x,\Delta^+(x)]}
\left| \int_{x}^{t}e^{\xi^2}\cos{e^{\xi^2}d\xi} \right| \quad \Longrightarrow
\end{equation*}
\begin{equation}\label{B:example:Phi}
\Phi(x) \le 2\max_{t\in [x,\Delta^+(x)]} \left| \int_{x}^{t}e^{\xi^2}\cos{e^{\xi^2}d\xi} \right|
\end{equation}
Now we find an estimate for integral in \eqref{B:example:Phi} for $|x|\gg 1$:
\begin{equation*}
\max_{t\in [x,\Delta^+(x)]} \left| \int_{x}^{t}e^{\xi^2}\cos{e^{\xi^2}d\xi} \right|
= \max_{t\in [x,\Delta^+(x)]} \left| \int_{x}^{t} \frac{d\sin{e^{\xi^2}}}{2\xi} \right|
\end{equation*}
\begin{equation*}
\le \max_{t\in [x,\Delta^+(x)]}
\left| \frac{\sin{e^{\xi^2}}}{2\xi}\Big|_{x}^{t} \right|
+ \frac{1}{2}\max_{t\in [x,\Delta^+(x)]}
\left| \int_{x}^{t}\frac{\sin{e^{\xi^2}}d\xi}{\xi^2} \right|
\end{equation*}
\begin{equation*}
\le \max_{t\in [x,\Delta^+(x)]} \left( \frac{1}{|t|}+ \frac{1}{|x|} \right)
\le \frac{c}{|x|}.
\end{equation*}
Finally,
$\Phi(x) \le \frac{c}{|x|} \textrm{ for } |x|\gg 1 \quad \Longrightarrow \quad \tilde{\varkappa}(x) \to 0 \textrm{ as } |x|\to \infty.$ The proof of \eqref{B:condition_e} for $t\in [\Delta^-,x]$ is analogous.

Now it remained to apply Theorem \ref{B:maintheorem} and finish a proof.
\qed


\section{Solution of Problem C}\label{C:proof} 
\setcounter{equation}{0}

\subsection{Proof of main result for the problem C}

In this section we prove Theorem \ref{C:maintheorem} (see \S\ref{C}).
We remind here that $D_p$ is the set of solutions $y\in L_p(R)$ of correctly solvable in $L_p(R)$ equation \eqref{theequation} with right side $f$ which belongs to the unit sphere $S_p = \{\,  f\in L_p : \|f\|_p = 1 \,\}$:
\begin{equation*}
D_p = \big\{\ y: -y'(x) + q(x)y(x) = f(x), \ y\in L_p(R),\ \|f\|_p=1\ \big\}
\end{equation*}
According to the proof of Theorem \ref{B:maintheorem} from \S\ref{B3:proof} the following equality holds:
\begin{equation}\label{C:asympt_ref}
G_p(x):=\sup_{y\in D_p^q}|y(x)| =
\left\{
\begin{array}{ll}
1\quad , & \quad p=1,\\
\displaystyle
\bigg( \int_{x}^{\infty} e^{-p' \int\limits_{x}^{t}q(\xi)d\xi} dt\bigg) ^{1/p'}, & \quad p\in(1, \infty),\\
\displaystyle 
\int_{x}^{\infty} e^{-\int\limits_{x}^{t}q(\xi)d\xi} dt\quad , & \quad p=\infty.
\end{array}
\right.
\end{equation}
\\
Now, our goal is a proof of Theorem \ref{C:maintheorem}. For convenience we give below its statement:
\\
{\bf Theorem \ref{C:maintheorem}.}
{\it
Let there exist $\alpha \ge 1$ and $\beta > 0$ such that for all $|x| \gg 1$ the inequality holds:
\begin{equation}\label{C:dt/dx-xgg1} 
\frac{1}{\alpha} \le \frac{d(t)}{d(x)} \le \alpha, \quad |t-x| \le \beta,
\end{equation}
where $d(x)$ is defined in \eqref{d}.\\
Then for $p\in (1,\infty]$ there exist $c, c(p) \in (0, \infty)$ such that for all $x\in R$ the inequalities hold:
\begin{equation}\label{C:mainresult} 
c^{-1}d(x)^{1/p'} \le G_p(x) \le c(p)d(x)^{1/p'}
\end{equation}
}
\paragraph{Proof.}
We remind here the following two facts:
\begin{enumerate}
\item[1)]
\begin{equation}\label{C:sup_d} 
d_0 < \infty, \quad d_0 = \sup_{x\in R} d(x) 
\end{equation}
Inequality \eqref{C:sup_d} follows from correct solvability of equation \eqref{theequation} in space $L_p(R),\ p\in [1,\infty]$ (see \S\ref{A:proof}).
\item[2)] The function $d(x)$ is positive and continuous for all $x\in R$.
\end{enumerate}
From 2) it follows that inequality \eqref{C:dt/dx-xgg1} holds for all $x\in R$ with substitution, possibly, of number $\alpha \ge 1$ for the larger $\alpha$.
Finally, the following inequality holds:
\begin{equation}\label{C:dt/dx-allx}
\frac{1}{\alpha} \le \frac{d(t)}{d(x)} \le \alpha, \quad |t-x| \le \beta, \quad x\in R.
\end{equation}
(For the new larger $\alpha$ we use the same symbol.)

Next, it is easy to see that from \eqref{C:sup_d} and \eqref{C:dt/dx-allx} it follows that
\begin{equation}\label{C:dt/dx-gamma}
\frac{1}{\gamma} \le \frac{d(t)}{d(x)} \le \gamma, \quad |t-x| \le d_0, \quad x\in R,
\end{equation}
where $\gamma(\ge \alpha \ge 1)$ is an absolute constant.

Let us now obtain the lower estimate in \eqref{C:mainresult}. Let $x\in R$ and
\begin{equation*}
f_x(t) =
\left\{
\begin{array}{ll}
\displaystyle
\frac{1}{d(x)^{1/p}}, & \quad t\in [x, x+d(x)], \\
\quad 0,& \quad t\notin [x, x+d(x)].
\end{array}
\right.
\end{equation*}
Then $\|f_x\|_p = 1$. Let $y_x$ be the solution of equation \eqref{theequation} ($y_x \in L_p(R)$) with the right side $f_x$. Then
\begin{eqnarray}\label{C:dt/dx-gamma-tag}
G_p(x) & \ge & |y_x(x)| = \left|\int_{x}^{\infty} f_x(t)
e^{-\int\limits_{x}^{t}q(\xi)d\xi}dt \right| \nonumber \\
&=& \int_{x}^{x+d(x)} \frac{1}{d(x)^{1/p}} \
e^{-\int\limits_{x}^{t}q(\xi)d\xi}dt \nonumber \\
&\ge& \frac{1}{d(x)^{1/p}} \
e^{-\int\limits_{x}^{x+d(x)}q(\xi)d\xi} \int_{x}^{x+d(x)} dt \nonumber \\
&\ge& d(x)^{-1/p} \ e^{-\int\limits_{x-d(x)}^{x+d(x)}q(\xi)d\xi}\ d(x)
= e^{-2} d(x)^{1/p'}.
\end{eqnarray}

We note that in order to obtain the lower estimate in \eqref{C:mainresult} we did not use any properties of $q$ except \eqref{q} and a requirement that fulfills a correct solvability of equation \eqref{theequation} in space $L_p(R),\ p\in [1,\infty]$.
The condition \eqref{C:dt/dx-xgg1} and the following from it condition \eqref{C:dt/dx-gamma} are used to obtain the upper estimate in \eqref{C:mainresult}.
Consider the integral:
\begin{equation*}
J_s(x)=\int_{x}^{\infty}e^{-s\int\limits_{x}^{t}q(\xi)d\xi}dt, \quad x\in R, \quad s\in (0, \infty)
\end{equation*}
We need the following definition:
\begin{definition}
We say that segments $\{\Delta_n\}_{n=1}^{\infty}$ constitute $R(x)$-covering of half-axis $[x, \infty)$ if they fulfill the following requirements:
\begin{enumerate}
\item[1)]
$\Delta_{n} = \big[x_n - d(x_n), x_n + d(x_n)\big] := \big[x_n - d_n, x_n + d_n\big] := \big[\Delta^{-}_{n}, \Delta^{+}_{n}\big], \ \ \ n\ge 1,$
\item[2)]
$\Delta_{n+1}^{-} = \Delta_{n}^{+}, \quad n\ge 1,$
\item[3)]
$\Delta_{1}^{-} = x,$
\item[4)]
$\displaystyle
\bigcup_{n\ge 1} \Delta_{n} = [x, \infty).$
\end{enumerate}
\end{definition}

\begin{lemma}
Let $d_0 < \infty$. Then for all $x\in R$ there exists $R(x)$-covering of half-axis $[x, \infty)$.
\end{lemma}
For the proof of the above lemma see \cite{CS4}.

\begin{lemma}
Let condition \eqref{C:dt/dx-xgg1} hold. Then
\begin{equation*}
J_s(x) \le c\, J_s^{*}(x), \quad x\in R,
\end{equation*}
where
\begin{equation*}
J_s^{*}(x)=\int_{x}^{\infty}e^{-\frac{s}{\gamma}\int\limits_{x}^{t}\frac{d\xi}{d(\xi)}}dt,
\end{equation*}
and $\gamma$ is defined in \eqref{C:dt/dx-gamma}.
\end{lemma}
\paragraph{Proof.}
\begin{eqnarray}\label{C:Jstar-le}
J_s(x) &=& \int_{x}^{\infty} e^{ -s \int\limits_{x}^{t} q(\xi)d\xi} dt =
\sum_{n=1}^{\infty} \int_{\Delta_n}
e^{ -s \int\limits_{x}^{t} q(\xi)d\xi }dt
\le 2 \sum_{n=1}^{\infty}
e^{ -s \int\limits_{\Delta_{1}^{-}}^{\Delta_{n}^{-}} q(\xi)d\xi }d_n \nonumber \\
&=& 2 \sum_{n=1}^{\infty} d_n e^{-2s(n-1)} \le 2e^2 \sum_{n=1}^{\infty} d_n e^{-2sn};
\end{eqnarray}
By \eqref{C:dt/dx-gamma} we obtain
\begin{eqnarray}\label{C:2n=sum-ge}
2n &=& \sum_{k=1}^{n} 2 = \sum_{k=1}^{n} \frac{1}{d_k} \int_{\Delta_k} d\xi =
 \sum_{k=1}^{n} \int_{\Delta_k} \frac{d(\xi)}{d_k} \frac{d\xi}{d(\xi)}\nonumber \\
&& \ge \frac{1}{\gamma} \sum_{k=1}^{n} \int_{\Delta_k} \frac{d\xi}{d(\xi)} = \frac{1}{\gamma} \int_{\Delta_{1}^{-}}^{\Delta_{n}^{+}} \frac{d\xi}{d(\xi)},
\end{eqnarray}
and by \eqref{C:2n=sum-ge} we continue to obtain the estimate in \eqref{C:Jstar-le}:
\begin{eqnarray*}
J_s(x) &\le& 2e^2 \sum_{n=1}^{\infty} d_n e^{-2sn}\le 2e^2 \sum_{n=1}^{\infty} d_n e^{-\frac{s}{\gamma}\int\limits_{\Delta_{1}^{-}}^{\Delta_{n}^{+}}\frac{d\xi}{d(\xi)}}
=e^2 \sum_{n=1}^{\infty} \int_{\Delta_n}e^{-\frac{s}{\gamma}
\int\limits_{\Delta_{1}^{-}}^{\Delta_{n}^{+}} \frac{d\xi}{d(\xi)}}dt\\
&=& e^2 \sum_{n=1}^{\infty} \int_{\Delta_n}
e^{-\frac{s}{\gamma}\int\limits_{\Delta_{1}^{-}}^{t}\frac{d\xi}{d(\xi)}}\cdot
e^{-\frac{s}{\gamma}\int\limits_{t}^{\Delta_{n}^{+}}\frac{d\xi}{d(\xi)}}dt
\le e^2 \sum_{n=1}^{\infty} \int_{\Delta_n}
e^{-\frac{s}{\gamma}\int\limits_{\Delta_{1}^{-}}^{t} \frac{d\xi}{d(\xi)}}dt \\
&=& e^2 \int_{x}^{\infty}
e^{-\frac{s}{\gamma}\int\limits_{x}^{t} \frac{d\xi}{d(\xi)}}dt = e^2 J_s^{*}(x).
\end{eqnarray*}

\begin{lemma}
Let condition \eqref{C:dt/dx-xgg1} hold. Then
\begin{equation}\label{C:J*_s-le}
J_s^{*} \le c \int_{x}^{x+d_0} e^{ -\frac{s}{\gamma} \int\limits_{x}^{t} \frac{d\xi}{d(\xi)}}dt, \quad x\in R.
\end{equation}
\end{lemma}
\paragraph{Proof.}
\begin{align}\label{C:J*_s=le}
& J_s^{*} =
\int_{x}^{\infty}  e^{ -\frac{s}{\gamma} \int\limits_{x}^{t} \frac{d\xi}{d(\xi)}}dt \nonumber \\
&=\int_{x}^{x+d_0}e^{ -\frac{s}{\gamma}\int\limits_{x}^{t}\frac{d\xi}{d(\xi)}}dt +
\int_{x+d_0}^{x+2d_0}e^{ -\frac{s}{\gamma}\int\limits_{x}^{t}\frac{d\xi}{d(\xi)}}dt +\ldots +
\int_{x+nd_0}^{x+(n+1)d_0}e^{ -\frac{s}{\gamma}\int\limits_{x}^{t}\frac{d\xi}{d(\xi)}}dt + \ldots \nonumber \\
&=\int_{x}^{x+d_0}e^{ -\frac{s}{\gamma}\int\limits_{x}^{t}\frac{d\xi}{d(\xi)}}dt 
\Bigg[ 1 + \sum_{n=1}^{\infty}
\int_{x+nd_0}^{x+(n+1)d_0}e^{ -\frac{s}{\gamma}\int\limits_{x}^{t}\frac{d\xi}{d(\xi)}}dt
\bigg( \int_{x}^{x+d_0}e^{ -\frac{s}{\gamma}\int\limits_{x}^{t}\frac{d\xi}{d(\xi)}}dt \bigg)^{-1} \Bigg] \nonumber \\
&\le\int_{x}^{x+d_0}e^{ -\frac{s}{\gamma}\int\limits_{x}^{t}\frac{d\xi}{d(\xi)}}dt
\left[ 1 + \sum_{n=1}^{\infty}
\frac{e^{ -\frac{s}{\gamma}\int\limits_{x}^{x+nd_0}\frac{d\xi}{d(\xi)}}\cdot d_0}
{e^{ -\frac{s}{\gamma}\int\limits_{x}^{x+d_0}\frac{d\xi}{d(\xi)}}\cdot d_0}
\right] \nonumber \\
&=\int_{x}^{x+d_0}e^{ -\frac{s}{\gamma}\int\limits_{x}^{t}\frac{d\xi}{d(\xi)}}dt
\left[ 1 + \sum_{n=1}^{\infty}
e^{ -\frac{s}{\gamma}
\left( \int\limits_{x}^{x+nd_0}\frac{d\xi}{d(\xi)} -
\int\limits_{x}^{x+d_0}\frac{d\xi}{d(\xi)}\right)}
\right]\nonumber \\
&= \int_{x}^{x+d_0}e^{ -\frac{s}{\gamma}\int\limits_{x}^{t}\frac{d\xi}{d(\xi)}}dt
\left[ 2 + \sum_{n=2}^{\infty}
e^{ -\frac{s}{\gamma}\int\limits_{x+d_0}^{x+nd_0}\frac{d\xi}{d(\xi)}} \right].
\end{align}
But from \eqref{C:dt/dx-gamma} and definition of $d_0$ it follows that
\begin{eqnarray}\label{C:int_do-nd0-lowest}
&&\int_{x+d_0}^{x+nd_0}\frac{d\xi}{d(\xi)}\nonumber \\
&=&\int_{x+d_0}^{x+2d_0}\frac{d\xi}{d(\xi)} + \int_{x+2d_0}^{x+3d_0}\frac{d\xi}{d(\xi)}
+ \ldots + \int_{x+(n-1)d_0}^{x+nd_0}\frac{d\xi}{d(\xi)} \nonumber \\
&=&\int_{x+d_0}^{x+2d_0}\frac{d(x+d_0)}{d(\xi)}\frac{d\xi}{d(x+d_0)} + \ldots \\
&+& \int_{x+(n-1)d_0}^{x+nd_0} \frac{d\big(x+(n-1)d_0\big)}{d(\xi)}
\frac{d\xi}{d\big(x+(n-1)d_0\big)} \nonumber \\
&\ge& \frac{1}{\gamma} \left[ \frac{d_0}{d(x+d_0)} + \frac{d_0}{d(x+2d_0)}
+\ldots + \frac{d_0}{d\big(x+(n-1)d_0\big)} \right]  \ge \frac{(n-1)}{\gamma}. \nonumber
\end{eqnarray}
Now \eqref{C:J*_s=le} and \eqref{C:int_do-nd0-lowest} yield:
\begin{equation*}
J_s^{*} \le \int_{x}^{x+d_0}
e^{ -\frac{s}{\gamma} \int\limits_{x}^{t} \frac{d\xi}{d(\xi)}} dt
\left[ 2 + \sum_{n=2}^{\infty} e^{ -\frac{s}{\gamma^2}(n-1)} \right]
= c \int_{x}^{x+d_0} e^{ -\frac{s}{\gamma} \int\limits_{x}^{t} \frac{d\xi}{d(\xi)}} dt.
\end{equation*}
\qed

\begin{lemma}
Let \eqref{C:dt/dx-xgg1} hold. Then
\begin{equation*}
J_s(x) \le c(s)d(x),\quad x\in R.
\end{equation*}
\end{lemma}
\paragraph{Proof.}
From \eqref{C:J*_s-le} and \eqref{C:dt/dx-xgg1} it follows that
\begin{eqnarray*}
& J_s^{*}(x) & \le
c \int_{x}^{x+d_0}e^{-\frac{s}{\gamma}\int\limits_{x}^{t}\frac{d\xi}{d(\xi)}}dt =
c \int_{x}^{x+d_0}e^{-\frac{s}{\gamma}\int\limits_{x}^{t}\frac{d(x)}{d(\xi)}\frac{d\xi}{d(x)}}dt \\
&& \le c \int_{x}^{x+d_0}e^{-\frac{s}{\gamma^2}\frac{t-x}{d(x)}}dt
\le c \int_{x}^{\infty}e^{-\frac{s}{\gamma^2}\frac{t-x}{d(x)}}dt
= \frac{c\gamma^2}{s}d(x).
\end{eqnarray*}
\qed
\\
Now we can easily obtain the upper estimate in \eqref{C:mainresult}:
\begin{equation*}
G_p^{\phantom{p}p'}(x)= \int_{x}^{\infty} e^{-p' \int\limits_{x}^{t}q(\xi)d(\xi)}dt
= J_{p'}(x) \le c J_{p'}^{*}(x) \le cd(x).
\end{equation*}
\qed

\subsection{Example}

\begin{example}
Consider the equation \eqref{theequation} with coefficient
\begin{equation}\label{C:example}
q(x) = e^{|x|} + e^{|x|}\cos{e^{2|x|}}, \quad x\in R.
\end{equation}
\end{example}
We need the following Lemma \ref{C:aux_lemma}:
\begin{lemma}{\rm\cite{CS4}}\label{C:aux_lemma}
Let condition \eqref{q} for $q$ hold. Suppose that there exist a continuous function $q_1(x)>0$ and a function $q_2\in L_1^{loc}(R)$ such that
\begin{equation*}
q(x) = q_1(x) + q_2(x), \quad x\in R,
\end{equation*}
and, besides,
\begin{equation}\label{C:aux_lemma:kappa12}
\sigma_1(x) \to 0,\, \sigma_2(x) \to 0, \quad |x|\to \infty,
\end{equation}
where
\begin{equation*}
\sigma_1(x) = \sup_{|z|\le 2/q_1(x)}
\left| \int_{0}^{z} \big[ q_1(x+t) - 2q_1(x) + q_1(x-t) \big] dt  \right|,
\end{equation*}
\begin{equation*}
\sigma_2(x) = \sup_{|z|\le 2/q_1(x)} \left| \int_{x-z}^{x+z}q_2(\xi)d\xi \right|.
\end{equation*}
Then
\begin{equation*}
q_1(x)d(x) = 1 + \varepsilon(x),\quad
|\varepsilon(x)| \le c \big[\sigma_1(x) + \sigma_2(x)\big], \ x\gg 1,
\end{equation*}
\begin{equation*}
c^{-1}\, \le\, q_1(x)d(x)\, \le\, c,\quad x\in R.
\end{equation*}
\end{lemma}
To apply Lemma \ref{C:aux_lemma} to \eqref{C:example} set
\begin{equation*}
q_1(x) = e^{|x|},\ q_2(x) = e^{|x|} \cos{e^{2|x|}}, \quad x\in R.
\end{equation*}
Since $q(x), q_1(x), q_2(x)$ are even, we consider only the case $x\ge 0$.
Let us find estimates for $\sigma_1(x)$ and $\sigma_2(x)$ ($x\gg 1$):
\begin{eqnarray*}
\sigma_1(x) &=& \sup_{|z|\le 2e^{-x}}
\left| \int_{0}^{z} \left( e^{x+t} - 2e^{x} + e^{x-t} \right) dt \right| \\
&=& e^{x} \sup_{|z|\le 2e^{-x}}
\left| \int_{0}^{z} \left( e^{t} - 2 + e^{-t} \right) dt \right| \\
&=& e^{x} \sup_{|z|\le 2e^{-x}}
\left| \int_{0}^{z} \left( t^2 + \frac{2t^4}{4!} + \ldots \right) dt \right| \\
&\le& c e^{x} \sup_{|z|\le 2e^{-x}} \int_{0}^{z} t^{2} dt = \frac{c}{e^{2x}};
\end{eqnarray*}
To find an estimate for $\sigma_2$ we use the second mean value theorem for integration:
\begin{eqnarray*}
\sigma_2(x) &=& \sup_{|z|\le 2e^{-x}}
\left| \int_{x-z}^{x+z} e^{t}\cos{e^{2t}} dt \right| \\
&=& \sup_{|z|\le 2e^{-x}}
\left| \int_{x-z}^{x+z} \frac{1}{2e^{t}}\Big[ 2e^{2t}\cos{e^{2t}}\Big]dt \right| \\
&\le& c \sup_{|z|\le 2e^{-x}} \frac {1}{e^{x}}
\left| \int_{x-z}^{\xi} 2e^{2t} \cos{e^{2t}}dt \right| = \frac{c}{e^{x}}.
\end{eqnarray*}
Therefore, since \eqref{C:aux_lemma:kappa12} holds, by Lemma \ref{C:aux_lemma} we obtain:
\begin{equation}\label{C:example:d(x)}
d(x) = \frac{1 + \varepsilon(x)}{e^{|x|}}, \quad |\varepsilon(x)| \le \frac{c}{e^{|x|}},\quad x\gg1.
\end{equation}
From \eqref{C:example:d(x)} it directly follows that $d_0 < \infty$ and inequalities \eqref{C:dt/dx-xgg1} hold. Then the equation \eqref{theequation} is correctly solvable in space $L_p(R), p\in [1, \infty]$ and inequalities \eqref{C:mainresult} are of the form:
\begin{equation*}
\frac{c^{-1}}{e^{|x|/p'}} \le G_p(x) \le \frac{c}{e^{|x|/p'}}, \quad p\in (1, \infty].
\end{equation*}


\section{Solution of Problem D}\label{D:proof}
\setcounter{equation}{0}

\subsection{Proof of main result for the problem of $\varepsilon$-strip}

In this section we prove Theorem \ref{D:maintheorem} (see \S\ref{D}). For convenience we remind below its statement.
\\ \\
{\bf Theorem \ref{D:maintheorem}.}
{\it
For $p=1$ the solutions of equation \eqref{theequation} do not tend in whole to zero as
$|x| \to \infty$. For $p\in (1, \infty]$ the solutions of equation \eqref{theequation} tend in whole to zero as $|x| \to \infty$ if and only if one of the following two equivalent conditions hold:
\begin{equation}\label{D:lim_int_q}
\lim_{|x| \to \infty} \int_{x-a}^{x+a} q(t)dt = \infty, \quad  \textrm{for all }\ a\in (0,\infty),
\end{equation}
\begin{equation}\label{D:lim_d}
\lim_{|x| \to \infty} d(x) = 0.
\end{equation}
}

\paragraph{Proof of Theorem \ref{D:maintheorem}.}
We need the following lemma.
\begin{lemma}\label{D:lemma-equiv_cond}
The conditions \eqref{D:lim_int_q} and \eqref{D:lim_d} are equivalent.
\end{lemma}
\paragraph{Proof.}
Let condition \eqref{D:lim_int_q} hold. Assume the contrary: $d(x)\nrightarrow 0, |x|\to \infty$.
Then there exist $\varepsilon > 0$ and a sequence $\{x_n\}_{n=1}^{\infty}$ such that
\begin{enumerate}
\item[1)] $|x_n| \to \infty, \quad n\to \infty$
\item[2)] $d(x_n)\ge \varepsilon, \quad n=1,2,\ldots$
\end{enumerate}
By 1) and 2) we obtain (see \eqref{D:lim_int_q}):
\begin{equation*}
2=\int_{x_n-d(x_n)}^{x_n+d(x_n)}q(t)dt \ge \int_{x_n-\varepsilon}^{x_n+\varepsilon}q(t)dt \to \infty,\quad n\to \infty,
\end{equation*}
contradiction, and therefore \eqref{D:lim_d} holds.\\
Now assume that \eqref{D:lim_d} holds. Fix $a>0$ and consider a segment $[x-a,x+a]$ as $|x|\to \infty$. Let $\{\Delta_n\}_{n=1}^{\infty}$ be a $R(x-a)$-covering of $[x-a, \infty)$. Since $(x-a)$ and $(x+a)$ tend to $\pm \infty$ as $x\to \pm\infty$, the lengths of $\Delta_n( n=1,2,\ldots,N_0)$ included in $[x-a,x+a]$ tend to zero and therefore $N_0\to \infty$ as $|x| \to \infty$. Consequently,
\begin{equation*}
\int_{x-a}^{x+a}q(t)dt\ge \sum_{k=1}^{N_0}\int_{\Delta_k}q(t)dt=\sum_{k=1}^{N_0}2 = 2N_0 \to \infty.
\end{equation*}
\qed
\begin{lemma}\label{D:lemma-Gp}
Let $p\in [1, \infty]$ and
\begin{equation*}
G_p(x) = \sup_{y\in D_p} |y(x)|, \quad x\in R.
\end{equation*}
Then the solutions of equation \eqref{theequation} tend in whole to zero as $|x|\to \infty$ if and only if
\begin{equation}\label{D:lim_Gp}
\lim_{|x|\to \infty} G_p(x) = 0.
\end{equation}
\end{lemma}
\paragraph{Proof.}
\textit{Necessity.}
Assume that solutions of equation \eqref{theequation} tend in whole to zero as $|x|\to \infty$. Then for any $\varepsilon > 0$ there exists $x_0(\varepsilon) \gg 1$ such that
\begin{equation*}
|y(x)| \le \varepsilon, \quad y\in D_p, \quad |x| \ge x_0(\varepsilon).
\end{equation*}
At the same time for $|x|\ge x_0(\varepsilon)$ we have
\begin{equation*}
G_p(x) = \sup_{y\in D_p} |y(x)| \le \varepsilon \quad \Longrightarrow \ \eqref{D:lim_Gp}
\end{equation*}
\textit{Sufficiency.} If condition \eqref{D:lim_Gp} holds then for any $\varepsilon > 0$ there exists $x_0(\varepsilon) \gg 1$ such that
\begin{equation*}
G_p(x) \le \varepsilon, \quad |x| \ge x_0(\varepsilon).
\end{equation*}
But then for any $y\in D_p$ we have
\begin{equation*}
y(x) \le \sup_{y\in D_p}|y(x)| = G_p(x) \le \varepsilon, \quad |x|\ge x_0(\varepsilon).
\end{equation*}
\qed
\begin{corollary}\label{D:Cor-G_1}
Let $p=1$. Then the solutions $y\in D_1$ do not tend in whole to zero as $|x|\to \infty$.
\end{corollary}
\paragraph{Proof.}
According to \eqref{C:asympt_ref}
\begin{equation*}
G_1(x) \equiv 1, \quad x\in R.
\end{equation*}
By Lemma \ref{D:lemma-Gp} we obtain the statement of Corollary \ref{D:Cor-G_1}. \qed
\paragraph{Proof of theorem \ref{D:maintheorem}.}\textit{Necessity.}
Let $p\in (1, \infty]$ and assume that solutions of equation \eqref{theequation} tend in whole to zero as $|x|\to \infty$. According to \eqref{C:dt/dx-gamma-tag} (see \S\ref{C:proof}) we have
\begin{equation*}
G_p(x) \ge e^{-2} d(x)^{1/p'}.
\end{equation*}
Then because of \eqref{D:lim_Gp} the condition \eqref{D:lim_d} holds and therefore (by Lemma \ref{D:lemma-equiv_cond}) the condition \eqref{D:lim_int_q} also holds.\\
\paragraph{Proof of theorem \ref{D:maintheorem}.}\textit{Sufficiency.}
Let \eqref{D:lim_d} hold. Let us show that \eqref{D:lim_Gp} also holds and prove our statement. We have (see \ref{C:asympt_ref}) for $\{\Delta_n\}_{n=1}^{\infty}$ - $R(x)$-covering of half-axis $[x,\infty)$:
\begin{eqnarray*}
G_p^{\phantom{p}p'}(x) &=&
\int_{x}^{\infty} e^{-p' \int\limits_{x}^{t}q(\xi)d\xi}dt \\
&=& \sum_{n=1}^{\infty} \int_{\Delta_n} e^{-p' \int\limits_{x}^{t}q(\xi)d\xi}dt\
\le \ \sum_{n=1}^{\infty}
\int_{\Delta_n} e^{-p' \int\limits_{\Delta_{1}^{-}}^{\Delta_{n}^{-}}q(\xi)d\xi}dt;
\end{eqnarray*}
We note that
\begin{eqnarray*}
\int_{\Delta_{1}^{-}}^{\Delta_{n}^{-}} q(\xi)d\xi \ = \
\sum_{k=1}^{n-1}\int_{\Delta_k} q(\xi)d\xi\ =\ \sum_{k=1}^{n-1}2 = 2(n-1), \quad n\ge 2\\
\Longrightarrow \quad \int_{\Delta_{1}^{-}}^{\Delta_{n}^{-}} q(\xi)d\xi\ \ge \ 2(n-1),\quad n\ge 1,
\end{eqnarray*}
and therefore,
\begin{equation*}
\sum_{n=1}^{\infty}
\int_{\Delta_n} e^{-p' \int\limits_{\Delta_{1}^{-}}^{\Delta_{n}^{-}}q(\xi)d\xi}dt\
\le \ \sum_{n=1}^{\infty} \frac{2d_n}{e^{2(n-1)p'}} =
c \sum_{n=1}^{\infty} \frac{d_n}{e^{2np'}}.
\end{equation*}
Denote
\begin{equation*}
S(x) = \sum_{n=1}^{\infty} \frac{d_n}{e^{2np'}}, \quad x\in R.
\end{equation*}
We consider cases $x\to +\infty$ and $x\to -\infty$ separately.\\
Let $x\to +\infty$. Since $d(x)\to 0$ as $|x|\to \infty$, for any $\varepsilon > 0$ there exists $x_0(\varepsilon)$ such that
\begin{equation*}
\sup_{x\ge x_0(\varepsilon)} d(x) \le \varepsilon.
\end{equation*}
Then for $x\ge x_0(\varepsilon)$ we obtain:
\begin{eqnarray*}
S(x) &=& \sum_{n=1}^{\infty} \frac{d_n}{e^{2np'}} \le
\sum_{n=1}^{\infty} \frac{\varepsilon}{e^{2np'}}=
\varepsilon \sum_{n=1}^{\infty} \frac{1}{e^{2np'}} = c(p)\varepsilon \\
&\Longrightarrow& \lim_{x\to \infty} S(x) = 0.
\end{eqnarray*}
\\
Let $x\to -\infty$. Fix $\varepsilon > 0$ and find $N(\varepsilon)$ such that
\begin{eqnarray*}
\sum_{n=N(\varepsilon)+1}^{\infty} \frac{d_n}{e^{2np'}} \le
d_0 \sum_{n=N(\varepsilon)+1}^{\infty} \frac{1}{e^{2np'}} \nonumber \\
= \frac{d_0}{e^{2(N(\varepsilon)+1)p'}} \sum_{k=0}^{\infty} \frac{1}{e^{2kp'}}
= \frac{c(p)d_0}{e^{2N(\varepsilon)p'}} \le \frac{\varepsilon}{2}.
\end{eqnarray*}
It is clear that for the above inequality to hold we have to choose $N(\varepsilon)$ such that ($d_0 = \sup\limits_{x\in R} d(x)$)
\begin{equation*}
N(\varepsilon) \ge \frac{1}{2p'} \ln \frac{2c(p)d_0}{\varepsilon}.
\end{equation*}
Enlarge, if needed, $x_0(\varepsilon)$ so that
\begin{equation*}
d(x) \le \frac{\varepsilon}{4N(\varepsilon)}, \quad |x| \ge x_0(\varepsilon).
\end{equation*}
Hence,
\begin{equation*}
d(x) \le \frac{\varepsilon}{4N(\varepsilon)}\quad \textrm { for }\ |x| \ge x_0(\varepsilon) + 2N(\varepsilon)d_0.
\end{equation*}
Segments $\Delta_n,\ n=1,\dots, N(\varepsilon)$ are included, obviously, in segment $[x, x+2N(\varepsilon)d_0]$ and therefore ($x_n$ is the center of $\Delta_n$, $x<0$):
\begin{equation*}
|x_n| \ge |x+2N(\varepsilon)d_0| \ge |x| - 2N(\varepsilon)d_0 \ge x_0(\varepsilon).
\end{equation*}
Thus, for $|x|\ge x_0(\varepsilon) + 2N(\varepsilon)d_0$ the inequality holds:
\begin{equation*}
d(x_n)\le \frac{\varepsilon}{4N(\varepsilon)},\quad n=1,2,\ldots,N(\varepsilon).
\end{equation*}
Then we have:
\begin{eqnarray*}
S(x) = \sum_{n=1}^{\infty} \frac{d_n}{e^{2np'}}
= \sum_{n=1}^{N(\varepsilon)} \frac{d_n}{e^{2np'}}
+ \sum_{N(\varepsilon) + 1}^{\infty} \frac{d_n}{e^{2np'}}
\le \sum_{n=1}^{N(\varepsilon)} \frac{\varepsilon}{4N(\varepsilon)} + \frac{\varepsilon}{2}
= \frac{3}{4}\varepsilon < \varepsilon
\end{eqnarray*}
\begin{equation*}
\Longrightarrow \quad \lim_{|x|\to \infty} S(x)=0 \quad
\Longrightarrow \quad \lim_{|x|\to \infty} G_p(x)=0.
\end{equation*}
We now apply Lemma \ref{D:lemma-Gp} to complete the proof. \qed
\begin{corollary}
Let condition \eqref{q} for $q$ hold. Suppose that there exist a continuous function $q_1(x)>0$ and a function $q_2\in L_1^{loc}(R)$ such that
\begin{equation*}
q(x) = q_1(x) + q_2(x), \quad x\in R,
\end{equation*}
and, besides,
\begin{equation*}
\sigma_1(x) \to 0,\, \sigma_2(x) \to 0, \quad |x|\to \infty,
\end{equation*}
where
\begin{equation*}
\sigma_1(x) = \sup_{|z|\le 2/q_1(x)}
\left| \int_{0}^{z} \big[ q_1(x+t) - 2q_1(x) + q_1(x-t) \big] dt  \right|,
\end{equation*}
\begin{equation*}
\sigma_2(x) = \sup_{|z|\le 2/q_1(x)} \left| \int_{x-z}^{x+z}q_2(t)dt \right|.
\end{equation*}
If at the same time $q_1(x)\to \infty$ as $|x|\to \infty$ then the solutions of equation \eqref{theequation} tend in whole to zero as $|x|\to \infty$.
\end{corollary}
\paragraph{Proof.}
The statement follows from Lemma \ref{C:aux_lemma} and Theorem \ref{D:maintheorem}.
\qed

\subsection{Examples}

\begin{example}
Consider the equation \eqref{theequation} with coefficient
\begin{equation}
q(x) = e^{|x|} + e^{|x|}\cos{e^{2|x|}}, \quad x\in R.
\end{equation}
\end{example}
As shown above (see \eqref{C:example:d(x)}), the following equality holds:
\begin{equation*}
d(x) = \frac{1 + \varepsilon(x)}{e^{|x|}}, \quad \lim_{|x|\to \infty}\varepsilon(x) = 0,
\end{equation*}
i.e. $d(x)\to 0$ as $|x|\to \infty$. Then by Theorem \ref{D:maintheorem} the solutions $y\in D_p$ of equation \eqref{theequation} in this case tend to zero in whole as $|x|\to \infty$.

\begin{example}
The equation \eqref{theequation} with coefficient
\begin{equation*}
q(x) = 1 + \cos(x), \quad x\in R,
\end{equation*}
has solutions that do not tend in whole to zero as $|x|\to \infty$.
\end{example}
\paragraph{Proof.}
Indeed, let
\begin{equation*}
x_k = (2k+1)\pi,\quad k=1,2,\ldots \quad \Longrightarrow
\end{equation*}
\begin{eqnarray*}
\int_{x_k - \frac{\pi}{2}}^{x_k + \frac{\pi}{2}}q(t)dt
&=& \pi + \int_{x_k - \frac{\pi}{2}}^{x_k + \frac{\pi}{2}} \cos(t)dt
= \pi + \sin(t)\big|^{x_k + \frac{\pi}{2}}_{x_k - \frac{\pi}{2}} = \pi -2 < 2 \\
&&\Longrightarrow \quad d(x_k) > \frac{\pi}{2},\quad k=1,2,\ldots \\
&&\Longrightarrow \quad d(x) \nrightarrow 0, \quad |x|\to \infty.
\end{eqnarray*}
By Theorem \ref{D:maintheorem} we obtain our statement. \qed


\section{Solution of Problem E}\label{E:proof}
\setcounter{equation}{0}

\subsection{Proof of main result for the problem E}

In this section we prove Theorem \ref{E:maintheorem} (see \S\ref{E}). For convenience we remind below its statement.
\\ \\
{\bf Theorem \ref{E:maintheorem}.}
{\it
Let condition \eqref{q} hold. Then for fixed $p\in[1,\infty]$ the operator $L^{-1}:L_p(R)\to L_p(R)$ is compact if and only if one of the following two equivalent conditions hold:
\begin{equation}\label{E:lim_int_q}
\lim_{|x| \to \infty} \int_{x-a}^{x+a} q(t)dt = \infty, \quad  \textrm{for all }\ a\in (0,\infty),
\end{equation}
\begin{equation}\label{E:lim_d}
\lim_{|x| \to \infty} d(x) = 0.
\end{equation}
}

\paragraph{Proof of Theorem \ref{E:maintheorem}.}
We remind that
\begin{equation*}
d(x)=\inf_{d>0}\ \Big\{d : \int_{x-d}^{x+d}q(t)dt=2\Big\},
\end{equation*}
and $d_0=\sup\limits_{x\in R}d(x)$. It is clear that $d_0 < \infty$. Let $\{\Delta_{n}\}_{n=1}^{\infty}$ be segments of $R(x)$-covering of half-axis $[x,\infty), x\in R$. Let $\lambda \in (0,\infty)$. Denote
\begin{equation*}
\begin{array}{c}
\displaystyle
I_1(x)=\int_{x}^{\infty}e^{-\int\limits_{x}^{t}(q(\xi)+\lambda)d\xi}dt,\quad x\in R,\\
\displaystyle
I_2(x)=\int^{x}_{-\infty}e^{-\int\limits^{x}_{t}(q(\xi)+\lambda)d\xi}dt,\quad x\in R.
\end{array}
\end{equation*}

\begin{lemma}\label{E:lemm-sup_JI_lambda}
For $\lambda \in (0,\infty)$ the following equalities hold:
\begin{equation}\label{E:sup_JI_lambda}
\begin{array}{c}
\displaystyle
\sup\limits_{x\in R} I_1(x) = \frac{1}{\lambda + \delta_1(\lambda)}, \quad
\delta_1{(\lambda)}>0, \\
\displaystyle
\sup\limits_{x\in R} I_2(x) = \frac{1}{\lambda + \delta_2(\lambda)}, \quad
\delta_2{(\lambda)}>0.
\end{array}
\end{equation}
\end{lemma}
\paragraph{Proof.}
Let us check the first equality in \eqref{E:sup_JI_lambda} (the second equality is proved in the same way). 
It is clear that (see conditions for correct solvability of equation \eqref{theequation}):
\begin{equation*}
1 = \int_{x}^{\infty} \big( q(t) + \lambda \big)
e^{-\int\limits_{x}^{t}(q(\xi) + \lambda) d\xi}dt, \quad x\in R.
\end{equation*}
Then
\begin{equation*}
1 = \int_{x}^{\infty} \big( q(t) + \lambda \big)
e^{-\int\limits_{x}^{t}(q(\xi) + \lambda) d\xi}dt
= \lambda I_1(x) +
\int_{x}^{\infty} q(t) e^{-\int\limits_{x}^{t}(q(\xi) + \lambda) d\xi}dt.
\end{equation*}
Now we find an upper estimate for the second term:
\begin{gather*}
\int_{x}^{\infty} q(t) e^{-\int\limits_{x}^{t}(q(\xi) + \lambda) d\xi}dt
= \sum_{n=1}^{\infty}
\int_{\Delta_n}q(t) e^{-\int\limits_{x}^{t}(q(\xi) + \lambda) d\xi}dt \\
\ge \sum_{n=1}^{\infty}
\int_{\Delta_n}q(t) e^{-\int\limits_{x}^{\Delta_{n}^{+}}(q(\xi) + \lambda) d\xi}dt
= \sum_{n=1}^{\infty}
e^{-\int\limits_{x}^{\Delta_{n}^{+}}(q(\xi)+\lambda)d\xi}
\cdot \int_{\Delta_n}q(t)dt; \\
\end{gather*}
We note here that
\begin{equation*}
\int_{\Delta_n}q(t)dt = \int_{x_n-d_n}^{x_n+d_n}q(t)dt = 2
= \frac{1}{d_n}\int_{x_n-d_n}^{x_n+d_n}dt = \frac{1}{d_n}\int_{\Delta_n}dt,
\end{equation*}
and then
\begin{gather*}
\sum_{n=1}^{\infty}
e^{-\int\limits_{x}^{\Delta_{n}^{+}}(q(\xi) + \lambda) d\xi}
\cdot \int_{\Delta_n}q(t)dt
= \sum_{n=1}^{\infty}
2e^{-\int\limits_{x}^{\Delta_{n}^{+}}(q(\xi) + \lambda) d\xi} \\
= \sum_{n=1}^{\infty}
e^{-\int\limits_{x}^{\Delta_{n}^{+}}(q(\xi) + \lambda) d\xi}
\frac{1}{d_n} \int_{\Delta_n} dt
\ge \frac{1}{d_0} \sum_{n=1}^{\infty}
\int_{\Delta_n} e^{-\int\limits_{x}^{\Delta_{n}^{+}}(q(\xi)+\lambda)d\xi} dt\\
= \frac{1}{d_0} \sum_{n=1}^{\infty}
\int_{\Delta_n} e^{-\int\limits_{x}^{t}(q(\xi)+\lambda)d\xi}
\cdot e^{-\int\limits_{t}^{\Delta_{n}^{+}}(q(\xi)+\lambda)d\xi} dt \\
\ge \frac{1}{d_0} \sum_{n=1}^{\infty}
\int_{\Delta_n} e^{-\int\limits_{x}^{t}(q(\xi)+\lambda)d\xi} dt
\cdot e^{-\int\limits_{\Delta_n}(q(\xi)+\lambda)d\xi}\\
= \frac{1}{d_0} \sum_{n=1}^{\infty}
\bigg( \int_{\Delta_n}e^{-\int\limits_{x}^{t}(q(\xi)+\lambda)d\xi} dt \bigg)
e^{-\int\limits_{\Delta_n}q(\xi)d\xi} \cdot e^{-2\lambda d_n} \\
\ge \frac{1}{d_0} \sum_{n=1}^{\infty}
\int_{\Delta_n}e^{-\int\limits_{x}^{t}(q(\xi)+\lambda)d\xi} dt
\cdot e^{-2-2\lambda d_0}
= \frac{e^{-2-2\lambda d_0}}{d_0} \sum_{n=1}^{\infty}
\int_{\Delta_n}e^{-\int\limits_{x}^{t}(q(\xi)+\lambda)d\xi} dt
\\
= \frac{e^{-2-2\lambda d_0}}{d_0}
\int_{x}^{\infty} e^{-\int\limits_{x}^{t}(q(\xi)+\lambda)d\xi} dt
= \frac{e^{-2-2\lambda d_0}}{d_0} I_1(x)
\end{gather*}
\begin{eqnarray*}
&\Longrightarrow& \quad 1 \ge \left( \lambda + \frac{e^{-2-2\lambda d_0}}{d_0} \right)  I_1(x) \\
&\Longrightarrow& \quad I_1(x) \ \le  \ \frac{1}{\lambda + e^{-2-2\lambda d_0}{d_0^{-1}}}, \quad x\in R \\
&\Longrightarrow& \quad \sup_{x\in R}  I_1(x) \ \le  \ \frac{1}{\lambda + e^{-2-2\lambda d_0}{d_0^{-1}}} \ < \ \frac{1}{\lambda}
\quad \Longrightarrow \quad \eqref{E:sup_JI_lambda}
\end{eqnarray*}
\qed

\begin{lemma}
For $p\in [1,\infty]$ the inequality holds:
\begin{equation}\label{E:||L+Elambda||^-1}
\|(L+\lambda E)^{-1}\|_{p\to p}\le \frac{1}{\lambda+\delta(\lambda)},\quad \delta(\lambda)>0,
\end{equation}
where $Ef \equiv f$ for any $f\in L_p(R)$.
\end{lemma}
\paragraph{Proof.}
Let $p\in (1, \infty)$. Since
\begin{equation*}
\left( (L + \lambda E)^{-1}f \right)(x) =
\int_{x}^{\infty} e^{-\int\limits_{x}^{t}(q(\xi) + \lambda) d\xi}f(t)dt, \quad x\in R, 
\end{equation*}
by H\"older inequality, Lemma \ref{E:lemm-sup_JI_lambda} and Fubini theorem we obtain:
\begin{gather*}
\|(L + \lambda E)^{-1}f\|_{p}^{p} = \int_{-\infty}^{\infty}
\bigg| \int_{x}^{\infty} e^{-\int\limits_{x}^{t}(q(\xi) + \lambda) d\xi}f(t)dt \bigg|^{p}dx\\
\le \int_{-\infty}^{\infty}
\bigg[ \int_{x}^{\infty} e^{-\int\limits_{x}^{t}(q(\xi) + \lambda) d\xi}|f(t)|dt \bigg]^{p}dx \\
= \int_{-\infty}^{\infty} \bigg[
\int_{x}^{\infty} e^{-\frac{1}{p'}\int\limits_{x}^{t} (q(\xi)+\lambda) d\xi}
\cdot e^{-\frac{1}{p}\int\limits_{x}^{t} (q(\xi)+\lambda) d\xi}|f(t)|dt
\bigg]^{p}dx \\
\le \int_{-\infty}^{\infty} \bigg[
\int_{x}^{\infty} e^{-\int\limits_{x}^{t} (q(\xi)+\lambda) d\xi}dt \bigg]^{p/p'}
\bigg[ \int_{x}^{\infty} e^{-\int\limits_{x}^{t} (q(\xi)+\lambda) d\xi}|f(t)|^{p} \bigg]dx\\
\le \sup_{x\in R} \bigg[
\int_{x}^{\infty} e^{-\int\limits_{x}^{t} (q(\xi)+\lambda) d\xi}dt \bigg]^{p/p'}
\cdot  \int_{-\infty}^{\infty} \bigg(
\int_{x}^{\infty} e^{-\int\limits_{x}^{t} (q(\xi)+\lambda) d\xi}|f(t)|^{p}dt \bigg)dx\\
= \frac{1}{(\lambda + \delta_1(\lambda))^{p/p'}}\int_{-\infty}^{\infty} \bigg(
\int_{x}^{\infty} e^{-\int\limits_{x}^{t} (q(\xi)+\lambda) d\xi}|f(t)|^{p}dt \bigg)dx\\
=  \frac{1}{(\lambda + \delta_1(\lambda))^{p/p'}} \int_{-\infty}^{\infty}|f(t)|^{p} \bigg(
\int_{-\infty}^{t} e^{-\int\limits_{x}^{t} (q(\xi)+\lambda) d\xi}dx \bigg)dt\\
\le \frac{1}{(\lambda + \delta_1(\lambda))^{p/p'}} \
\sup_{t\in R} \bigg( \int^{t}_{-\infty} e^{-\int\limits_{x}^{t} (q(\xi)+\lambda) d\xi}dx \bigg)
\cdot \int_{-\infty}^{\infty} |f(t)|^{p}dt \\
= \frac{1}{(\lambda + \delta_1(\lambda))^{p/p'}} \
\frac{1}{(\lambda + \delta_2(\lambda))} \ \|f\|^{p}_{p}
\end{gather*}
\begin{eqnarray*}
&\Longrightarrow& \quad \|(L + \lambda E)^{-1}f\|_p \le
\frac{1}{(\lambda + \delta_1(\lambda))^{1/p'}} \
\frac{1}{(\lambda + \delta_2(\lambda))^{1/p}} \ \|f\|_{p} \\
&\Longrightarrow& \quad \eqref{E:||L+Elambda||^-1}
\end{eqnarray*}

Let $p=1$ or $p=\infty$. Then (see \cite[ch.V, \S2, 4-5]{KA})
\begin{eqnarray*}
\|(L + \lambda E)^{-1}\|_{1\to 1} &=& \sup_{t\in R}
\int^{t}_{-\infty} e^{-\int\limits_{x}^{t} (q(\xi)+\lambda) d\xi}dx \\
\|(L + \lambda E)^{-1}\|_{C(R)\to C(R)} &=& \sup_{x\in R}
\int_{x}^{\infty} e^{-\int\limits_{x}^{t} (q(\xi)+\lambda) d\xi}dt
\end{eqnarray*}
In these cases \eqref{E:||L+Elambda||^-1} follows immediately from Lemma \ref{E:lemm-sup_JI_lambda}. \qed

\begin{lemma}\label{E:Lemma=L^-1=ab}
For $\lambda\in (0,\infty)$ the operator equality holds:
\begin{equation*}
L^{-1} = (L + \lambda E)^{-1} (E + S_{\lambda})^{-1},
\end{equation*}
where
\begin{equation*}
S_{\lambda} : L_p(R) \to L_p(R), \quad \|S_{\lambda}\|_{p\to p} < 1.
\end{equation*}
\end{lemma}
\paragraph{Proof.}
Let $f\in L_p(R)$ and
\begin{equation*}
Ly =  f,\quad y\in L_p(R).
\end{equation*}
Then
\begin{equation*}
(L + \lambda E) y = f + \lambda y.
\end{equation*}
Therefore, $f + \lambda y \in L_p(R)$ (since $y\in L_p(R)$),
\begin{eqnarray*}
\Longrightarrow& \quad y = (L + \lambda E)^{-1} f + \lambda(L + \lambda E)^{-1} y, \\
\Longrightarrow& \quad y - \lambda (L + \lambda E)^{-1} y =  (L + \lambda E)^{-1} f, \\
\Longrightarrow& \quad [E - \lambda (L + \lambda E)^{-1}] y = (L + \lambda E)^{-1} f.
\end{eqnarray*}
Denote
\begin{equation*}
S_{\lambda} = -\lambda (L + \lambda E)^{-1}, \quad \lambda > 0.
\end{equation*}
By \eqref{E:||L+Elambda||^-1} we have:
\begin{equation*}
\|S_{\lambda}\|_{p\to p} = \lambda \|(L + \lambda E)^{-1}\|_{p\to p}
\le \frac{\lambda}{\lambda + \delta(\lambda)} < 1.
\end{equation*}
Hence,
\begin{equation*}
(E + S_{\lambda})y = (L + \lambda E)^{-1}f,
\end{equation*}
and $\|S_{\lambda}\|_{p\to p} < 1$. Therefore, the operator $(E + S_{\lambda})^{-1}$ exists and is bounded (see \cite{LS}). Thus we obtain that for all $f\in L_p(R),\ p\in [1,\infty]$ the following relations hold:
\begin{eqnarray*}
\left\{
\begin{array}{ll}
y = (E + S_{\lambda})^{-1} (L + \lambda E)^{-1} f, \\
y = L^{-1} f.
\end{array}
\right.
\end{eqnarray*}
\begin{eqnarray}
&\Longrightarrow& \quad L^{-1} f = (E + S_{\lambda})^{-1} (L + \lambda E)^{-1} f, \quad
\forall f\in L_p(R), \nonumber \\
&\Longrightarrow& \quad L^{-1} = (E + S_{\lambda})^{-1} (L + \lambda E)^{-1}.
\end{eqnarray}
Since $  (E + S_{\lambda})^{-1} = \sum\limits_{n=0}^{\infty} (-S_{\lambda})^{n}, \quad S_{\lambda} = -\lambda (L + \lambda E)^{-1}$, the operators $S_{\lambda}$ and $(L + \lambda E)^{-1}$ are invertible. And we obtain
\begin{equation}\label{E:L^-1=ab=ba}
L^{-1} = (E + S_{\lambda})^{-1} (L + \lambda E)^{-1} =
(L + \lambda E)^{-1} (E + S_{\lambda})^{-1}.
\end{equation}
\qed
\\
Let us now refer to the proof of Theorem \ref{E:maintheorem}.
\paragraph{Proof of Theorem \ref{E:maintheorem}.}\textit{Necessity.}
Let $\lambda = 1$ and let operator $L^{-1}: L_p(R) \to L_p(R)$ be compact. From \eqref{E:L^-1=ab=ba} it follows that
\begin{equation*}
(L + \lambda E)^{-1} = L^{-1} (E + S_{\lambda})
\end{equation*}
Since $(E + S_{\lambda})$ is bounded, $(L + \lambda E)^{-1} : L_p(R) \to L_p(R)$ is compact. Then by Theorem \ref{E:prevtheorem} condition \eqref{E:lim_int_q} holds.\\
\paragraph{Proof of Theorem \ref{E:maintheorem}.}
\textit{Sufficiency.}
Let \eqref{E:lim_int_q} hold. Then by Theorem \ref{E:prevtheorem} the operator $(L + \lambda E)^{-1}$ is compact. From \eqref{E:L^-1=ab=ba} it follows that
\begin{equation*}
L^{-1} = (L + \lambda E)^{-1} (E + S_{\lambda})^{-1}.
\end{equation*}
Since $(E + S_{\lambda})^{-1}$ is bounded, by Theorem \ref{E:prevtheorem} we obtain that operator $L^{-1}$ is compact. \qed


\end{document}